\documentclass[11pt]{article}
\usepackage{amssymb,amsmath,latexsym, color}

\allowdisplaybreaks

\newtheorem{thm}{Theorem}[section]

\newtheorem{lemma}[thm]{Lemma}
\newtheorem{prop}[thm]{Proposition}

\newtheorem{remark}[thm]{Remark}
\numberwithin{equation}{section}

\def\pf{{\medskip\noindent {\bf Proof. }}}
\def\qed{{\hfill $\Box$ \bigskip}}

\def\sA {{\cal A}}

 \def\bH {{\mathbb H}}

  \def\bR {{\mathbb R}}

\def\nn{\nonumber}
\def\R {{\mathbb R}}
\def\wt{\widetilde}
\def\wh{\widehat}

\def\E{{\mathbb E}}
\def\P{{\mathbb P}}

\def\bea{\begin{align*}}
\def\eea{\end{align*}}
\def\bee{\begin{equation}}
\def\eee{\end{equation}}

\def\eps{\varepsilon}

\def\wh{\widehat}

\def\1{{\bf 1}}

\makeatletter
\@addtoreset{equation}{section}

\makeatother

\begin{document}
\bibliographystyle{plain}

\title{\Large \bf Global Heat Kernel Estimates for $\Delta+\Delta^{\alpha/2}$ \\
in Half-space-like Domains }

\author{
{\bf Zhen-Qing Chen}\thanks{Research partially supported
by NSF Grants    DMS-0906743 and DMR-1035196.},
\quad
{\bf Panki Kim}\thanks{Research supported by Basic Science Research Program through the National
Research Foundation of Korea(NRF) grant funded by the Korea
government(MEST)(2010-0001984).} \quad and  \quad {\bf Renming Song}
}
\date{(February 23, 2011)}

\maketitle

\begin{abstract}
Suppose that $d\ge 1$ and $\alpha\in (0, 2)$. In this paper, by
using probabilistic methods, we establish sharp two-sided pointwise
estimates for the Dirichlet heat kernels of $\{\Delta+ a^\alpha
\Delta^{\alpha/2}; \ a\in (0, 1]\}$ on half-space-like $C^{1, 1}$
domains for all time $t>0$. The large time estimates for
half-space-like domains are very different from those for bounded
domains. Our estimates are uniform in $a \in (0, 1]$ in the sense
that the constants in the estimates are independent of $a\in (0,
1]$. Thus it yields the Dirichlet heat kernel estimates for Brownian
motion in half-space-like domains by taking $a\to 0$. Integrating
the heat kernel estimates with respect to the time variable $t$, we obtain
uniform sharp two-sided estimates for the Green functions
 of $\{\Delta+ a^\alpha \Delta^{\alpha/2}; \ a\in (0, 1]\}$ in
half-space-like $C^{1, 1}$ domains in $\bR^d$.
\end{abstract}



\bigskip
\section{Introduction and Setup}

This paper is a natural continuation of \cite{CKS5} where small time
sharp two-sided estimates for the Dirichlet heat kernel of $\Delta+
\Delta^{\alpha/2}$ on any  $C^{1, 1}$ open sets and large time
sharp two-sided estimates for bounded $C^{1,1}$ open sets are
obtained. In this paper we give sharp two-sided estimates for the
Dirichlet heat kernel of $\Delta+ \Delta^{\alpha/2}$ on
half-space-like $C^{1, 1}$ domains for all time. The large time
Dirichlet heat kernel estimates for half-space-like domains are very
different from those for bounded open sets. See below for the
definition of half-space-like $C^{1, 1}$ open sets.

Throughout this paper, we assume that $d\ge 1$ is an integer and
$\alpha\in (0, 2)$. Let $X^0=(X^0_t,\, t\ge 0)$ be a Brownian motion
in $\bR^d$ with generator $\Delta=\sum_{i=1}^d
\frac{\partial^2}{\partial x_i^2}$, and $Y=(Y_t,\, t\ge 0)$ be an
independent (rotationally) symmetric $\alpha$-stable process in
$\bR^d$ whose generator is the fractional Laplacian $\Delta^{\alpha
/2}$. For $u \in C^{\infty}_c(\bR^d)$, the space of smooth functions
with compact support, the fractional Laplacian can be written in the
form
\begin{equation}\label{e:1.1}
\Delta^{\alpha /2} u(x)\, =\, \lim_{\eps \downarrow 0}\int_{\{y\in
\bR^d: \, |y-x|>\eps\}} (u(y)-u(x)) \frac{\sA (d, \alpha)
}{|x-y|^{d+\alpha}}\, dy,
\end{equation}
where $ {\cal A}(d, \alpha):= \alpha2^{\alpha-1}\pi^{-d/2}
\Gamma(\frac{d+\alpha}2) \Gamma(1-\frac{\alpha}2)^{-1}. $ Here
$\Gamma$ is the Gamma function defined by $\Gamma(\lambda):=
\int^{\infty}_0 t^{\lambda-1} e^{-t}dt$ for every $\lambda > 0$.

For any $a>0$, we define $X^a$ by $X_t^a:=X^0_t+ a Y_t$. We will
call the process $X^a$ the independent sum of the Brownian motion
$X^0$ and the symmetric $\alpha$-stable process  $Y$ with weight
$a>0$.
The L\'evy process $X^a$ is uniquely
determined by its characteristic function
$$
\E_x \left[ e^{i\xi\cdot(X^a_t-X^a_0)}
\right]\,=\,e^{-t(|\xi|^{2}+a^{\alpha}|\xi|^{\alpha})} \qquad
\hbox{for every } x\in \bR^d \hbox{ and }  \xi\in \bR^d
$$
and
its infinitesimal generator
 is $\Delta+ a^\alpha \Delta^{\alpha/2}$.
Since
$$
a^{\alpha}|\xi|^{\alpha}= \int_{\bR^d}(1-\cos(\xi\cdot y))\, \frac{
a^{\alpha} {\cal A}(d, \alpha)}{|y|^{d+\alpha}}dy,
$$
$X^a$ has L\'evy intensity function
$$
J^a(x, y)=j^a(|x-y|) :=a^{\alpha}{\cal A}(d,
\alpha)|x-y|^{-(d+\alpha)}.
$$
The  function  $J^a (x, y)$  determines a L\'evy system for $X^a$,
which describes the jumps of the process $X^a$: for any non-negative
measurable function $f$ on $\bR_+ \times \bR^d\times \bR^d$ with
$f(s, y, y)=0$ for all $y\in \bR^d$, any stopping time $T$ (with
respect to the filtration of $X^a$) and any $x\in \bR^d$,
\begin{equation}\label{e:levy}
\E_x \left[\sum_{s\le T} f(s,X^a_{s-}, X^a_s) \right]= \E_x \left[
\int_0^T \left( \int_{\bR^d} f(s,X^a_s, y) J^a(X^a_s,y) dy \right)
ds \right]
\end{equation}
(see, for example, \cite[Proof of Lemma 4.7]{CK} and \cite[Appendix
A]{CK2}).

Let $p^a(t, x, y)$ be the transition density of
$X^a$ with respect to the Lebesgue measure on $\bR^d$.
The function $p^a(t, x, y)$ is smooth on $(0, \infty)\times
\bR^d\times \bR^d$. For any $\lambda>0$,
$(\lambda X^{a}_{\lambda^{-2} t}, t\geq 0)$ has the same
distribution as $(X^{a \lambda^{(\alpha-2)/\alpha}}_t, t\ge 0)$ (see
the second paragraph of
\cite[Section 2]{CKS5}), so we
have
\begin{equation}\label{scalingrd}
p^{a\lambda^{(\alpha-2)/\alpha}} ( t,  x, y)= \lambda^{-d} p^{a}
(\lambda^{-2}t, \lambda^{-1} x, \lambda^{-1} y) \qquad \hbox{for }
t>0 \hbox{ and } x, y \in \bR^d.
\end{equation}

For $a>0$ and $C>0$, define
\begin{equation}\label{eq:qd}
h^a_C(t, x, y):
 =  \left(t^{-d/2} \wedge (a^\alpha t)^{-d/\alpha}
\right) \wedge \left(t^{-d/2} e^{-C |x-y|^2/t}+\left( (a^{\alpha}
t)^{-d/\alpha}\, \wedge \frac{a^{\alpha}t}{|x-y|^{d+\alpha}} \right)
\right).
\end{equation}
Here and in the sequel, we use ``$:=$" as a way of definition and,
for $a, b\in \bR$, $a\wedge b:=\min \{a, b\}$ and $a\vee b:=\max\{a,
b\}$. The following sharp two-sided estimates on $p^a(t, x, y)$
follow from \eqref{scalingrd} and the main results in \cite{CK08,
SV07} that give the sharp estimates on $p^1(t, x, y)$.

\begin{thm}\label{T:1.1}
There are constants $c, C_1\geq 1$ such that, for all $a\in [0,
\infty)$ and $(t, x, y)\in (0, \infty]\times \bR^d\times \bR^d$
$$
c^{-1} \, h^a_{C_1}(t,x, y)\leq p^a (t, x, y) \leq c\,
h^a_{1/C_1}(t, x, y).
$$
\end{thm}

We record a simple but  useful observation. Its proof will be
given at the end of this section.

\begin{prop}\label{C:1.2}
For every $c>0$ and $c_1>0$, there is a constant $c_2\geq 1$ such
that for any $a>0$,
$$
c_2^{-1} \left( (a^{\alpha} t)^{-d/\alpha}\, \wedge
\frac{a^{\alpha}t}{|x-y|^{d+\alpha}} \right) \leq h^a_{c}(t, x, y)
\leq c_2  \left( (a^{\alpha} t)^{-d/\alpha}\, \wedge
\frac{a^{\alpha}t}{|x-y|^{d+\alpha}} \right)
$$
holds when either $t\geq c_1 a^{-2\alpha/(2-\alpha)}$ or $|x-y|\geq
a^{-\alpha/(2-\alpha)}$.
\end{prop}

Recall that a domain (an connected open set) $D$ in $\bR^d$ (when
$d\ge 2$) is said to be $C^{1,1}$ if there exist a localization
radius $ R_0>0 $ and a constant $\Lambda_0>0$ such that for every
$z\in\partial D$, there exist a $C^{1,1}$ function $\psi=\psi_z:
\bR^{d-1}\to \bR$ satisfying $\psi(0)=0$, $ \nabla\psi (0)=(0,
\dots, 0)$, $\| \nabla \psi \|_\infty \leq \Lambda_0$, $| \nabla
\psi (x)-\nabla \psi (z)| \leq \Lambda_0 |x-z|$, and an orthonormal
coordinate system $CS_z$: $y=(y_1, \cdots, y_{d-1}, y_d):=(\wt y, \,
y_d)$ with origin at $z$ such that $B(z, R_0 )\cap D= \{y=({\tilde
y}, y_d) \in B(0, R_0) \mbox{ in } CS_z: y_d > \psi (\wt y) \}$. The
pair $( R_0, \Lambda_0)$ will be called the $C^{1,1}$
characteristics of
the domain $D$.

For an open set $D\subset \bR^d$ and $x\in D$, we will use
$\delta_D(x)$ to denote the Euclidean distance between $x$ and
 $D^c$. For a domain  $D\subset \bR^d$ and $\lambda_0 \geq 1$,
 we say {\it the path distance in $D$ is comparable to the Euclidean
distance with characteristic $\lambda_0$} if   for every $x, y\in
D$, there is a rectifiable curve $l$ in $D$ connecting $x$ to $y$ so
that the length of $l$ is no larger than $\lambda_0|x-y|$. Clearly,
such a property holds for all bounded $C^{1,1}$ domains, $C^{1,1}$
domains with compact complements and domains above the graphs of
bounded $C^{1,1}$ functions.

For any open subset $D\subset \bR^d$, we use $\tau^a_D$ to denote
the first time the process $X^a$ exits $D$. We define the process
$X^{a,D}$ by $X^{a,D}_t=X^a_t$ for $t<\tau^a_D$ and
$X^{a,D}_t=\partial$ for $t\ge \tau^a_D$, where $\partial$ is a
cemetery point. $X^{a,D}$ is called
the subprocess of $X^a$ in $D$.
 The generator of $X^{a,D}$ is $(\Delta+ a^\alpha
\Delta^{\alpha/2})|_D$. It follows from \cite{CK08} that $X^{a,D}$
has a continuous transition density $p^a_D(t, x, y)$ with respect to
the Lebesgue measure.

One can easily see that, when $D$ is bounded, the operator
$-(\Delta + a^\alpha \Delta^{\alpha/2})|_D$ has discrete spectrum.
In this case, we use $\lambda^{a, D}_1>0$ to denote the smallest
eigenvalue of $-(\Delta + a^\alpha \Delta^{\alpha/2})|_D$.

The following is a particular case of a more general result proved
in \cite[Theorem 1.3]{CKS5} (cf. Proposition \ref{C:1.2} above).

\begin{thm}\label{t:main}
Suppose that $D$ is a $C^{1,1}$ domain in $\bR^d$ with
characteristics $( R_0, \Lambda_0)$ such that the path distance in
$D$ is comparable to the Euclidean distance with characteristic
$\lambda_0$.
\begin{description}
\item{\rm (i)}
For every $M>0$ and $T>0$, there are constants $c_1=c_1( R_0,
\Lambda_0, \lambda_0, M, \alpha, T)\geq 1$ and $C_2=C_2( R_0,
\Lambda_0,  \lambda_0, M, \alpha, T)\geq 1$ such that for all $a
\in (0, M]$ and $(t, x, y)\in (0, T]\times D\times D$,
\begin{eqnarray*}
 && \hskip -0.6truein
c_1^{-1}\,\left(1\wedge \frac{\delta_D(x)}{\sqrt{t}}
\right)\left(1\wedge \frac{\delta_D(y)}{\sqrt{t}}  \right)
h^a_{C_2}(t, x, y) \\
&& \le  p^a_D(t, x, y) \leq c_1 \left(1\wedge
\frac{\delta_D(x)}{\sqrt{t}} \right)\left(1\wedge
\frac{\delta_D(y)}{\sqrt{t}}  \right)h^a_{1/C_2}(t, x, y).
\end{eqnarray*}

\item{\rm (ii)}
Suppose in addition that $D$ is bounded. For every
$M>0$ and $T>0$, there is a constant $c_2=c_2(D, M, \alpha, T)\geq
1$ so that for all $a \in (0, M]$ and $(t, x, y)\in [T,
\infty)\times D\times D$,
$$
c_2^{-1}\, e^{- t \, \lambda^{a, D}_1 }\, \delta_D (x)\, \delta_D
(y) \,\leq\,  p^a_D(t, x, y) \,\leq\, c_2\, e^{- t\, \lambda^{a,
D}_1}\, \delta_D (x) \,\delta_D (y).
$$
\end{description}
\end{thm}

Note that Theorem \ref{t:main} does not give large time estimates
for $p^a_D(t, x, y)$ when $D$ is unbounded. The goal of this paper
is to establish two-sided large time estimates on $p^a_D(t, x, y)$
for a large class of unbounded $C^{1, 1}$ domains, namely
half-space-like $C^{1, 1}$ domains. A domain $D$ is said to be
half-space-like if, after isometry, there exist two real numbers
$b_1\le b_2$ such that $\bH_{b_2}\subset D\subset \bH_{b_1}$. Here
and throughout this paper, $\bH_b$ stands for the set $\{x=(x_1,
\dots, x_d)\in \R^d: x_d>b\}$. We will denote $\bH_0$ by $\bH$.

Now we are in a position to state the main result of this paper.
For $a>0$, define $\phi_a(r):= r\wedge (r/a)^{\alpha/2}$.

\begin{thm} \label{t:main_hsl}
Suppose $D$ is a half-space-like $C^{1,1}$ domain with $C^{1, 1}$
characteristic $(R_0, \Lambda_0)$ and $\bH_b\subset D\subset \bH$
for some $b>0$ such that the path distance in $D$ is comparable to
the Euclidean distance with characteristic $\lambda_0$.
 Then for any $M\geq 1$, there exist constants
$c_i=c_i
 ( R_0, \Lambda_0,  \lambda_0, M, \alpha, b)
\geq 1$, $i=1, 2,$ such that for all $a\in (0, M]$ and $(t, x,y) \in
(0, \infty) \times D\times D$,
\begin{align}
&c_1^{-1} \left(1\wedge \frac{\phi_a(\delta_{D}(x))}{\sqrt{t}}\right)
\left( 1\wedge \frac{\phi_a(\delta_{D}(y))}{\sqrt{t}}\right)
 h^a_{c_2}  (t, x, y) \nn
\\
\le&   p^a_D(t, x, y) \leq c_1 \left(1\wedge \frac{\phi_a(\delta_{D}(x))}{\sqrt{t}}\right)
\left( 1\wedge \frac{\phi_a(\delta_{D}(y))}{\sqrt{t}}\right)
 h^a_{1/c_2}(t, x, y).
\label{e:newform}
\end{align}
\end{thm}

\begin{remark}\label{R:1.5} \rm
(i) The L\'evy exponent for $X^a$ is $\Phi_a(|\xi|)$
with
 $\Phi_a (r):=r^2+ a^\alpha r^\alpha$. The function $\phi_a(r)$ is
related to $\Phi_a(r)$ as follows.
$$
\frac1{\Phi_a(1/r)}=\frac{1}{r^{-2}+a^\alpha r^{-\alpha}} \asymp
\frac1{r^{-2}} \wedge \frac1{(a/r)^\alpha} =  r^2 \wedge
(r/a)^{\alpha}  = \phi_a(r)^2.
$$
Here for two non-negative functions $f$ and $g$, the notation $f\asymp g$
means that there is a  positive constant $c\geq 1$ so that
$g( x)/c \leq f (x)\leq c g(x)$ in the common domain of definition
for $f$ and $g$.
Hence in view of Theorem \ref{T:1.1}, the estimate \eqref{e:newform}
can be restated as follows. For every $M>0$, there are
constants $c_1, c_2\geq 1$ so that for every $a\in (0, M]$ and
$(t, x, y)\in (0, \infty)\times D \times D$,
\begin{eqnarray}
&& c_1^{-1} \left( 1\wedge \frac{1}{ t \Phi_a (1/\delta_D(x))} \right)^{1/2}
\left( 1\wedge \frac{1}{ t \Phi_a (1/\delta_D(y))} \right)^{1/2}
 p^a (t, c_2 x, c_2 y) \nn
\\
&\le&   p^a_D(t, x, y) \leq c_1 \left( 1\wedge \frac{1}{ t \Phi_a (1/\delta_D(x))} \right)^{1/2} \left( 1\wedge \frac{1}{ t \Phi_a (1/\delta_D(y))} \right)^{1/2}
 p^a (t, x/c_2, y/c_2).
\end{eqnarray}
We conjecture that the above Dirichlet heat kernel estimates
hold for a large class of rotationally symmetric L\'evy processes in $\bR^d$;
see \cite[Conjecture]{CKS4}.

\medskip

\noindent (ii)
 Note that  $t\leq a^{2\alpha/(\alpha-2)}$ if and only
if   $(a^\alpha t)^{-d/\alpha} \geq t^{-d/2}$.
If $(\delta_D(x)/a)^{\alpha/2}<\delta_D(x)$, then $\delta_D (x)\geq
a^{\alpha/(\alpha -2)}$ and so $\delta_D(x)\wedge
(\delta_D(x)/a)^{\alpha/2} \geq a^{\alpha/(\alpha -2)}$. Thus when
$t\leq a^{2\alpha/(\alpha-2)}$ and
$(\delta_D(x)/a)^{\alpha/2}<\delta_D(x)$, we have $\frac{
(\delta_D(x)/a)^{\alpha/2}}{\sqrt{t}} \geq \frac{ a^{\alpha/(\alpha
-2)}}{a^{\alpha/(\alpha-2)}}=1$, and consequently
$$ 1\wedge \frac{\delta_D(x) \wedge (\delta_D(x)/a)^{\alpha/2}}{\sqrt{t}}
=1=1\wedge \frac{\delta_D(x)}{\sqrt{t}}.
$$
Hence in view of Theorem \ref{T:1.1} and
Proposition \ref{C:1.2},
the statement of Theorem \ref{t:main_hsl} can be restated as follows.
For all $a\in (0, M]$ and
$(t, x, y)\in (0,  a^{2\alpha/(\alpha-2)}]\times D\times D$,
\begin{align}
& c_1^{-1} \left(1\wedge \frac{\delta_D(x)}{\sqrt{t}}\right)
\left(1\wedge \frac{\delta_D(y)}{\sqrt{t}}\right) \left( t^{-d/2}
e^{-c_2 |x-y|^2/t}+ t^{-d/2} \wedge \left( \frac{a^\alpha
t}{|x-y|^{d+\alpha}}\right)\right)\nn\\
&\le p^a_D(t, x, y) \leq c_1 \left(1\wedge
\frac{\delta_D(x)}{\sqrt{t}} \right)\left(1\wedge
\frac{\delta_D(y)}{\sqrt{t}}  \right) \left( t^{-d/2}
e^{-|x-y|^2/(c_2t)}+ t^{-d/2} \wedge \left( \frac{a^\alpha
t}{|x-y|^{d+\alpha}}\right)\right) \label{e:old_1}
\end{align}
and for all $a\in (0, M]$ and $(t, x, y)\in [
a^{2\alpha/(\alpha-2)}, \infty)\times D\times D$,
\begin{align}
&c_1^{-1} \left(1\wedge \frac{\delta_{D}(x) \wedge(a^{-1}\delta_D
(x))^{\alpha/2}}{\sqrt{t}}\right)\left(1\wedge \frac{\delta_{D}(y)
\wedge(a^{-1}\delta_D (y))^{\alpha/2}}{\sqrt{t}}\right)
\left((a^\alpha t)^{-d/\alpha}\wedge\frac{a^\alpha
t}{|x-y|^{d+\alpha}} \right)\nn\\
&   \le   p^a_D(t, x, y) \leq \nn\\
&c_1 \left(1\wedge \frac{\delta_{D}(x) \wedge(a^{-1} \delta_D
(x))^{\alpha/2}}{\sqrt{t}}\right) \left(1\wedge \frac{\delta_{D}(y)
\wedge(a^{-1} \delta_D (y))^{\alpha/2}}{\sqrt{t}}\right)
\left((a^\alpha t)^{-d/\alpha}\wedge\frac{a^\alpha
t}{|x-y|^{d+\alpha}} \right).\label{e:old_2}
\end{align}
In fact, Theorem  \ref{t:main_hsl} will be proved in this form.
\qed
\end{remark}

\begin{remark}\label{R:1.6}
\rm
Unlike \cite{CKS4, CT}, there are dramatic differences between the
behavior of the heat kernel $p^a_D(x, y)$ on half-space-like
 $C^{1,1}$ domains and disconnected half-space-like $C^{1,1}$ open sets
even if $x$ and $y$ are in the same connected component. For
example, if $D$ is $\bH \cup B(x_0,1)$ where $x_0=(0, \dots, 0, -2)$
and $x,y \in B(x_0,1)$, then, as $a \to 0$, $p^a_D(x, y)$ converges
to $p^0_{B(x_0,1)}(x, y)$, the Dirichlet heat kernel for Brownian
motion on $B(x_0,1)$. Thus, in this case, the heat kernel estimates
for
 $p^a_D(t, x, y)$ when $t$ is large  cannot
be of the form \eqref{e:newform} even if
$x$ and $y$ are in the same connected component. Furthermore, as one
can see from \cite[Theorem 1.3]{CKS5}, when $D$ is a disconnected
half-space-like $C^{1,1}$ open set (containing bounded connected
component), we can not expect that the heat kernel estimates for
$p^a_D(x, y)$ to be written in a simple form as the one in
\eqref{e:newform}. To keep our exposition as transparent as
possible, we are content with establishing the heat kernel estimates
for half-space-like $C^{1,1}$ domains. \qed
\end{remark}

Integrating the heat kernel estimates in Theorem \ref{t:main_hsl}
with respect to $t$, we get sharp two-sided estimates on the Green
function $G^a_D(x, y):=\int_0^\infty p^a_D(t, x, y)dt$
for $X^a$ in  half-space-like $C^{1, 1}$ domains $D$.

Define for $d\geq 1$ and $a>0$,
\begin{equation}\label{e:f^a}
f^a_D(x, y) = \begin{cases} \frac{1} {|x-y|^{d-\alpha}}
\left(a^{-\alpha/2} \wedge \frac{  \phi_a (\delta_D(x))}{
|x-y|^{\alpha/2}   }\right) \left( a^{-\alpha/2} \wedge \frac{
\phi_a (\delta_D(y))}{ |x-y|^{\alpha/2}}
\right)  &\hbox{when } d>\alpha, \smallskip  \\
\log \left( \left( 1+ a \, \frac{ \phi_a(\delta_D(x))  \phi_a(
\delta_D(y)) }{ |x-y| }\right)^{1/a}\right)
&\hbox{when } d=1=\alpha,  \smallskip \\
\frac{ \phi_a(\delta_D(x)) \phi_a ( \delta_D(y)) }{ |x-y|  } \wedge
\left(  a^{-1}  \left( \phi_a (\delta_D(x))  \phi_a (
\delta_D(y))\right) ^{(\alpha-1)/\alpha} \right) &\hbox{when }
d=1<\alpha .
\end{cases}
\end{equation}
For $d\ge 2$ and $a>0$, define
$$
g_D^a (x, y) = \begin{cases} \frac{1} {|x-y|^{d-2}} \left(1\wedge
\frac{ \delta_D(x) \delta_D(y)}{ |x-y|^{2}}\right)
\quad &\hbox{when } d\geq 3, \smallskip \\
 \log\left(1+\frac{a^{2\alpha/(\alpha-2)}\wedge (\delta_D(x)\delta_D(y))}{ |x-y|^{2}}\right)
\quad &\hbox{when } d=2,
\end{cases}
$$
for $d=1$ and  $a>0$, define
$$
g_D^a (x, y)=
\begin{cases} \left( \delta_{D}({x}) \delta_{D} ({y})\right)^{1/2}
   \wedge \frac{  \delta_{D}({x}) \delta_{D} ({y})}{ |{x}-{y}| }
   \wedge  \left(a^{-\alpha}( \delta_{D}({x}) \delta_{D} ({y}))^{(\alpha-1)/2} \right)
&\hbox{when } \alpha \in (1, 2)  , \smallskip  \\
\frac{
\delta_{D}({x}) \delta_{D} ({y})}{ |{x}-{y}| } \wedge
  \log\left(1+a\left( \delta_{D}({x})\delta_{D}({y})\right)^{1/2} \right)^{1/a}
&\hbox{when } \alpha=1,\smallskip \\
\left( \delta_{D}({x}) \delta_{D} ({y})\right)^{1/2} \wedge \frac{
\delta_{D}({x}) \delta_{D} ({y})}{ |{x}-{y}| } \wedge a^{\alpha/(\alpha-2)}
&\hbox{when } \alpha \in (0,1).
\end{cases}
$$

\begin{thm}\label{t:gf-estimates}
Suppose $D$ is a half-space-like $C^{1,1}$ domain with $C^{1, 1}$
characteristic $(R_0, \Lambda_0)$ and $\bH_b\subset D\subset \bH$
for some $b>0$ such that the path distance in $D$ is comparable to
the Euclidean distance with characteristic $\lambda_0$. Then for any
$M>0$, there exists a constant $c=c(M, R_0, \Lambda_0, \lambda_0, b,
\alpha)\geq 1$ such that for all $a\in (0, M]$ and $(x,y) \in
D\times D$,
\begin{eqnarray}
&&c^{-1}g^a_D(x,y)\le G^a_D(x,y) \le c g^a_D(x,y)\qquad \text{when }
|x-y| \le a^{-\alpha/(2-\alpha)},
\label{e:gf-estimates1}\\
&&c^{-1}f_D^a(x,y)\le G^a_D(x,y) \le  c f_D^a(x,y) \qquad \text{when }
|x-y|\ge a^{-\alpha/(2-\alpha)}\, .
 \label{e:gf-estimates2}
\end{eqnarray}
\end{thm}

\begin{remark}\label{R:1.9} \rm
(i) Note that, when $d \ge 3$, $g^a_D(x, y)$ is independent of $a$ and
is comparable to the Green function of
Brownian motion in a bounded $C^{1,1}$ domain or
in a domain above the graph of a bounded $C^{1,1}$ function.
On the other hand, when $d \le 2$, $g^a_D(x, y)$ depends on $a$, which
is due to recurrent nature of one- and two-dimensional Brownian motion.

\medskip
\noindent (ii)
 Observe that if $(X^{a, D}_t, t\geq 0)$ is the  subprocess in $D$ of
the independent sum of a Brownian motion and a symmetric
$\alpha$-stable process in $\R^d$ with weight $a$, then $(\lambda
X^{a, D}_{\lambda^{-2} t}, t\geq 0)$ is the subprocess in $\lambda
D$ of the independent sum of a Brownian motion and a symmetric
$\alpha$-stable process in $\R^d$ with  weight $a
\lambda^{(\alpha-2)/\alpha}$ (see the second paragraph of
\cite[Section 2]{CKS5}). Consequently for any $\lambda>0$, we have
\begin{equation}\label{e:scaling}
p^{a\lambda^{(\alpha-2)/\alpha}}_{\lambda D} ( t,  x, y)=
\lambda^{-d} p^{a}_D (\lambda^{-2}t, \lambda^{-1} x, \lambda^{-1} y) \qquad
\hbox{for } t>0 \hbox{ and } x, y \in
\lambda D.
\end{equation}
When $D$ is a half space, we see from \eqref{e:scaling} that
 Theorems \ref{t:main_hsl} and \ref{t:gf-estimates} hold
 with $M=\infty$.

\medskip \noindent (iii)
 The estimates in Theorems \ref{t:main_hsl} and \ref{t:gf-estimates}
 are uniform in
 $a\in (0, M]$
 in the sense
that the constants $c_1$, $c_2$ and $c$ in the estimates
 are independent of $a\in (0, M]$. Since $X^a$ converges weakly to
 $X^0$, by taking $a\to 0$
 these estimates yield the following estimates for the heat kernel
 $p^0_D(t, x, y)$ and Green function $G^0(x, y)$ of
 Brownian motion in half-space-like domains $D$ in which
the path distance is comparable to the Euclidean distance:
\begin{align}
\hskip -1truein & c_1^{-1} \left(1\wedge
\frac{\delta_D(x)}{\sqrt{t}}\right) \left(1\wedge
\frac{\delta_D(y)}{\sqrt{t}}\right) \, t^{-d/2} e^{-c_2
|x-y|^2/t} \nonumber \\
 &\le  p^0_D(t, x, y)\leq c_1 \left(1\wedge
\frac{\delta_D(x)}{\sqrt{t}} \right)\left(1\wedge
\frac{\delta_D(y)}{\sqrt{t}}  \right) \, t^{-d/2}
e^{-|x-y|^2/(c_2t)}  \label{e:1.12}
\end{align}
for every $(t, x, y)\in (0, \infty)\times D \times D$,
and
\begin{equation}\label{e:1.13}
 c_2^{-1}\,
g^0_D(x,y)\,\le\, G^0_D(x,y)\, \le\,
\,c_2\, g^0_D(x,y)\qquad
\hbox{for } x, y \in D.
\end{equation}
The estimates \eqref{e:1.12} and \eqref{e:1.13} extend the main
results in \cite{So}, where the corresponding estimates were
established for domains in $\R^d$ with $d\geq 3$ that are above the
graphs of  bounded $C^{1,1}$ functions.

\medskip

\noindent
 (iv)  By Theorem \ref{t:main_hsl}, the boundary decay rate of the
 Dirichlet heat kernel
of $\Delta + \Delta^{\alpha/2}$ is  given by
$1\wedge \frac{\delta_{D}(x) \wedge \delta_D (x)^{\alpha/2}}{\sqrt{t}}$.
This indicates that the Dirichlet heat kernel estimates
for $\Delta + \Delta^{\alpha/2}$ in half-space-like $C^{1,1}$ domains
cannot be   obtained by a ``simple" perturbation argument from $\Delta$
nor from  $\Delta^{\alpha/2}$.
\end{remark}

The main difficulty of this paper is to obtain the correct boundary
decay rate of the Dirichlet heat kernel of $\Delta +
\Delta^{\alpha/2}$.
In \cite{CKS5}, the correct boundary decay rate for small $t$
was established by using some exit distribution estimates obtained in
\cite{CKSV}. Unfortunately the estimates in \cite{CKSV} are not
suitable for the present case.
Thus, in this paper we give some
different forms of exit distribution estimates that are suitable for
large time estimates.  The first step is, similar to \cite{BBC, G,
CKSV}, to compute $(\Delta +  \Delta^{\alpha/2}) h$ for certain test
functions. But unlike \cite{CKSV}, we do not use combinations of
test functions to serve as subharmonic and superharmonic functions
to obtain our desired estimates.
Instead, we use a generalization of Dynkin's  formula
 to obtain the desired exit distribution estimates directly. We
believe that our approach to obtain the correct boundary decay rate
is quite general and  may be used for other types of jump processes.

Throughout this paper, the constants $C_1, C_2, C_3$, $R_0, R_1,
R_2, R_3$ will be fixed. The lower case constants $c_1, c_2, \dots$
will denote generic constants whose exact values are
 not important and can change from one appearance to another.
The dependence of the lower case constants on the dimension $d$ will
not be mentioned explicitly. We will use $\partial$ to denote a
cemetery point and for every function $f$, we extend its definition
to $\partial$ by setting $f(\partial )=0$. We will use $dx$ or
$m(dx)$ to denote the Lebesgue measure in $\bR^d$. For a Borel set
$A\subset \bR^d$, we also use $|A|$ to denote its $d$-dimensional
Lebesgue measure. For every function $f$, let $f^+:=f \vee 0$.

In the remainder of this paper we will always assume that $D$ is a
half-space-like $C^{1,1}$ domain with $C^{1, 1}$ characteristic
$(R_0, \Lambda_0)$ and $\bH_b\subset D\subset \bH$ for some $b>0$
such that the path distance in $D$ is comparable to
the Euclidean distance with characteristic $\lambda_0$
and that $t_0$,
$x_0$ and $y_0$ are described as below.

Fix $t_0 \ge b^2$ and let $e_d$ be the unit vector in the direction
of the $x_d$-axis. For $x$ and $y$ in $D$, define the points
\begin{align} \label{e:x0y0}
x_0 := x + 2t_0^{1/2}e_d \qquad\text{and}\qquad y_0 := y +
2t_0^{1/2}e_d \ .
\end{align}
Observe that
\begin{equation}\label{e:geH1}
\delta_D(x_0) \ge
\delta_{\bH}
(x_0) > t_0^{1/2}, \qquad
\delta_D(y_0) \ge  \delta_{\bH}
(y_0) > t_0^{1/2},
\end{equation}
and
$ |x-x_0| = |y-y_0| =2 t_0^{1/2}$.
Note that when $D=\bH$, we can take $t_0$ to be any positive number.
Now as a consequence of Theorem \ref{t:main}, we have the following
result.

\begin{lemma}\label{L:ratio}
There  exists $c = c (b, t_0, R_0, \Lambda_0, \alpha,
\lambda_0)\geq 1$ such that for all
$x, z\in D$,
\begin{align}
c^{-1} \left( 1\wedge \delta_D(x) \right) \le
\frac{p^1_D(t_0,x,z)}{p^1_D(t_0,x_0,z)} \le c \left( 1\wedge
\delta_D(x) \right) \label{E:ratio-half-space-like1}.
\end{align}
\end{lemma}

\pf
Let $C_2$ be the constant in Theorem \ref{t:main} (i) with $T=t_0$.
From Proposition \ref{C:1.2} and Theorem \ref{t:main} (i), it is easy to see that
\begin{equation}\label{e:sv07}
h^1_{C_2}(t_0, x, y)\asymp 1\wedge \frac{1}{|x-y|^{d+\alpha}}\quad
\mbox{ and } \quad h^1_{1/C_2}(t_0, x, y)\asymp 1\wedge
\frac{1}{|x-y|^{d+\alpha}}.
\end{equation}
By Theorem \ref{t:main} (i) and \eqref{e:geH1}, we see that
\begin{align}
c_1^{-1}\left( 1\wedge \frac{\delta_D(x)}{\sqrt{t_0}} \right) \left(
\frac{h^1_{C_2}(t_0, x, z) }{h^1_{1/C_2}(t_0, x_0, z)} \right) \le
\frac{p^1_D(t_0,x,z)}{p^1_D(t_0,x_0,z)} \le c_1\left( 1\wedge
\frac{\delta_D(x)}{\sqrt{t_0}} \right) \left( \frac{h^1_{1/C_2}(t_0,
x, z) }{h^1_{C_2}(t_0, x_0, z)} \right). \label{e:eded1}
\end{align}

For $z \in B(x_0,2^{-1} t_0^{1/2})$ we have
$$
\frac32  t_0^{1/2} \le |x_0-x|- |z-x_0| \le |x-z| \le  |z-x_0|
+|x_0-x| = |z-x_0| + 2 t_0^{1/2}\le \frac 52  t_0^{1/2} .
$$
Similarly, for $z \in B(x,2^{-1} t_0^{1/2})$ we have $ \frac32
t_0^{1/2} \le |x-z_0| \le  \frac52  t_0^{1/2}.$ Thus in these cases,
\eqref{E:ratio-half-space-like1} follows from \eqref{e:eded1}.

In the case $z \not\in B(x,2^{-1} t_0^{1/2})\cup B(x_0,2^{-1}
t_0^{1/2})$,  we have $ |x-z| \le  |z-x_0| +|x_0-x| = |z-x_0| + 2
t_0^{1/2} \le 5|z-x_0| $ and $ |x_0-z| \le  |z-x| +|x_0-x| = |z-x| +
2 t_0^{1/2} \le 5|z-x|. $ So $5^{-1} |x_0-z| \le |z-x|\le 5
|x_0-z|$. Therefore by \eqref{e:sv07}
$$
\frac{h^1_{1/C_2}(t_0, x, z) }{h^1_{C_2}(t_0, x_0, z)}  \le c_2
\quad \text{and} \quad \frac{h^1_{C_2}(t_0, x, z) }{h^1_{1/C_2}(t_0,
x_0, z)}  \ge c_3.
$$
\qed

\begin{lemma}\label{L:rationew}
For any $M>0$, there exists $c =c (b, t_0, R_0, \Lambda_0, \alpha,
 \lambda_0) \geq 1$ such that for all  $a\in (0, M]$ and $x, z\in
D$,
\begin{eqnarray}
&&\hskip -0.6truein   c^{-1} \left( 1\wedge
 \delta_D(x)\right)\left( 1\wedge
\delta_D(z) \right)h^a_{25C_2}(t_0, x_0, z)\nonumber\\
&&\le p^a_D(t_0,x,z) \le c  \left( 1\wedge \delta_D(x) \right)\left(
1\wedge \delta_D(z) \right)h^a_{1/(25C_2)}(t_0, x_0, z)
\label{E:ratio-half-space-likenew1}
\end{eqnarray}
where  $C_2$ is the constant in Theorem \ref{t:main} (i) with $T=t_0$.
\end{lemma}

\pf
By Theorem \ref{t:main} (i), we see that
 \bee\label{e:eded4}
c_1^{-1} \left( 1\wedge \delta_D(x) \right)\left( 1\wedge
\delta_D(z) \right)h^a_{C_2}(t_0, x, z)\le p^a_D(t_0,x,z) \le c_1
\left( 1\wedge \delta_D(x) \right)\left( 1\wedge \delta_D(z)
\right)h^a_{1/(C_2)}(t_0, x, z).
 \eee
By the same argument as in the proof of Lemma \ref{L:ratio},
 $\frac32  t_0^{1/2} \le |x-z| \le  \frac 52  t_0^{1/2}$ for $z \in
B(x_0,2^{-1} t_0^{1/2})$,  $ \frac32  t_0^{1/2} \le |x-z_0| \le
\frac52  t_0^{1/2}$ for $z \in B(x,2^{-1} t_0^{1/2})$,  and $5^{-1}
|x_0-z| \le |z-x|\le 5 |x_0-z|$ for $z \not\in B(x,2^{-1}
t_0^{1/2})\cup B(x_0,2^{-1} t_0^{1/2}).$ The assertion of the lemma
follows by considering each cases in \eqref{e:eded4}. \qed

The following elementary result will play an important role later in
this paper. Recall that $D$, $t_0$, $x_0$ and $y_0$ are described as
above.

\begin{lemma}\label{l:piw2}
For any $t_0 \ge b^2$ and $M>0$, there exists a constant
$c=c(\alpha, M, t_0, b)>1$
 such that for any $a\in (0, M]$ and $(t,
 x)\in [t_0, \infty)\times D$,
\begin{eqnarray*}
(1\wedge\delta_D(x))\left(1\wedge\frac{\delta_{\bH}(x_0)
\wedge(a^{-1}\delta_{\bH}(x_0))^{\alpha/2}}{\sqrt{t}}\right)&\le&
c\left(1\wedge \frac{\delta_{D}(x)
\wedge(a^{-1}\delta_{D}(x))^{\alpha/2}}{\sqrt{t}}\right),\\
(1\wedge\delta_D(x))\left(1\wedge\frac{\delta_{\bH_b}(x_0)
\wedge(a^{-1}\delta_{\bH_b}(x_0))^{\alpha/2}}{\sqrt{t}}\right)&\ge&
c^{-1}\left(1\wedge \frac{\delta_{D}(x)
\wedge(a^{-1}\delta_{D}(x))^{\alpha/2}}{\sqrt{t}}\right).
\end{eqnarray*}
\end{lemma}

\pf Note that
\begin{align*}
\delta_D(x) +   t_0^{1/2} \le \delta_{\bH_b}(x_0) \le \delta_D(x) +
2t_0^{1/2} \quad\text{and}\quad \delta_D(x)+2t_0^{1/2} \le
\delta_{\bH}(x_0) \le \delta_D(x) + 3t_0^{1/2}.
\end{align*}
When $\delta_D(x) > t_0^{1/2}$, we have
$
\delta_D(x) \le \delta_{\bH_b}(x_0) < \delta_{\bH}(x_0) \le
4\delta_D(x).
$
Thus in this case, the conclusion of the lemma is trivial. From now
on, we assume that $\delta_D(x)\le t^{1/2}_0$. In this case,
using the fact
$t \ge t_0$
 and $a\in (0, M]$, we have
\begin{eqnarray*}
&&(1\wedge\delta_D(x))\left(1\wedge\frac{\delta_{\bH}(x_0)
\wedge(a^{-1}\delta_{\bH}(x_0))^{\alpha/2}}{\sqrt{t}}\right)
\asymp (1\wedge\delta_D(x))\left(1\wedge\frac{\delta_{\bH_b}(x_0)
\wedge(a^{-1}\delta_{\bH_b}(x_0))^{\alpha/2}}{\sqrt{t}}\right)\\
&&\asymp \delta_D(x)\left(1\wedge \frac1{\sqrt{t}}\right)\asymp
1\wedge \frac{\delta_D(x)}{\sqrt{t}}\asymp 1\wedge \frac{\delta_{D}(x)
\wedge(a^{-1}\delta_{D}(x))^{\alpha/2}}{\sqrt{t}}.
\end{eqnarray*}
The proof is now complete. \qed

\bigskip

\noindent{\bf Proof of
Proposition \ref{C:1.2}.} We first deal with
the case $a=1$. For $t\geq c_1$ and $r\ge 0$,
$$
t^{-d/2} e^{-
c r^2/t} \leq t^{-d/2} \frac{c_2}{(c r^2/t)^{(d+\alpha)/2}}
\leq c_3 \frac{t^{\alpha/2}}{r^{d+\alpha}}\leq c_4 \frac{t}{r^{d+\alpha}}.
$$
Hence for $t\geq c_1$,
$$   t^{-d/\alpha} \wedge \left( t^{-d/2}e^{-
cr^2/t}
+ t^{-d/\alpha}\, \wedge
\frac{t}{|x-y|^{d+\alpha}} \right) \asymp t^{-d/\alpha}\, \wedge
\frac{t}{|x-y|^{d+\alpha}}.
$$
Thus $h^1_c (t, x, y)\asymp t^{-d/\alpha}\, \wedge
\frac{t}{|x-y|^{d+\alpha}}$ on $[c_1, \infty)\times \R^d\times
\R^d$. On the other hand, for  $r\geq 1$,
$$
t^{-d/2} e^{-
cr^2/t} \leq t^{-d/2} \frac{c_5}{(
cr^2/t)^{(d/2)+1}}
= \frac{c_6t}{r^{d+2}} \leq \frac{c_6t}{r^{d+\alpha}}.
$$
So for $t\in (0, c_1]$ and $r\geq 1$,
$$ t^{-d/2}e^{- cr^2/t}
+ \left( t^{-d/2}\wedge \frac{t}{r^{d+\alpha}}\right)
\asymp  t^{-d/2} \wedge \frac{t}{r^{d+\alpha}}
\asymp \frac{t}{r^{d+\alpha}}\asymp   t^{-d/\alpha} \wedge \frac{t}{r^{d+\alpha}}.
$$
Thus we conclude that $h^1_c (t, x, y) \asymp t^{-d/\alpha}\, \wedge
\frac{t}{|x-y|^{d+\alpha}}$ for $t\leq c_1$ and $|x-y|\geq 1$. In
summary, we have
\begin{equation}\label{e:dssd}
h^1_c (t, x, y)\asymp t^{-d/\alpha}\, \wedge
\frac{t}{|x-y|^{d+\alpha}}
\end{equation}
when $t\geq c_1$ or $|x-y|\geq 1$. For $a>0$,
with $\lambda=a^{\alpha/(2-\alpha)}$, by \eqref{e:dssd}
\begin{eqnarray*}
h^a_c (t, x, y)&=& \lambda^d h^1_c(\lambda^2t, \lambda x, \lambda y)  \\
 &\asymp& \lambda^d \left( \left(\lambda^2 t\right)^{-d/\alpha} \wedge
 \frac{\lambda^2 t}{\lambda^{d+\alpha} |x-y|^{d+\alpha}}  \right)
 = (a^\alpha t)^{-d/\alpha} \wedge \frac{a^\alpha t}{|x-y|^{d+\alpha}},
 \end{eqnarray*}
provided either $\lambda^2t\geq c_1$ or $\lambda |x-y|\geq 1$.
This proves the proposition.
\qed

\section{Preliminary estimates}\label{S:int}

We will focus on the case $D=\bH$ in Sections
\ref{S:int}--\ref{sec:lower-half-space}. In this section we will
prove some preliminary estimates that will be used to establish our
heat kernel estimates in $\bH$. We start with some
 one-dimensional results.

Let $S$ be the sum of a unit drift and an $\alpha/2$-stable
subordinator and let $W$ be an independent one-dimensional Brownian
motion. Define a process $Z$ by $Z_t=W_{S_t}$. The process $Z$ is
simply the process $X^1$ in the case of dimension 1 defined in the
previous section. We will use the fact that $S$ is a complete
subordinator, that is, the L\'evy measure of $S$ has a completely
monotone density (for more details see
\cite{SSV} or \cite{SV06}).
Let $\overline{Z}_t:=\sup\{0\vee Z_s:0\le s\le t\}$ and let $L_t$ be
a local time of $\overline{Z}-Z$ at $0$. $L$ is also called a local
time of the process $Z$ reflected at the supremum. Then the right
continuous inverse $L^{-1}_t$ of $L$ is a
 subordinator and is called the ladder time process of $Z$. The
process $\overline{Z}_{L^{-1}_t}$ is also a
 subordinator and is called the ladder height process of $Z$.
(For the basic properties of the ladder time and ladder height
processes, we refer our readers to \cite[Chapter 6]{Ber}.) Let
$V(dr)$ denote the potential measure of the ladder height process
$\overline{Z}_{L^{-1}_t}$ of $Z$ and $v(r)$ its density, which is a
decreasing function on $[0, \infty)$. We know by \cite[(5.1)]{KSV09}
that
\begin{equation}\label{e:estimates-for-v}
v(r) \asymp 1 \wedge
r^{\alpha/2-1} \quad \text{for } r > 0.
\end{equation}
Let $G_{(0,\infty)}$ be the Green function of $Z^{(0,\infty)}$, the
subprocess of $Z$ in $(0,\infty)$.
 By using \cite[Theorem 20, p.~176]{Ber} which was originally proved
in \cite{Sil}, the following formula for $G_{(0,\infty)}$  was shown
in \cite[Proposition 2.8]{KSV}:
\begin{equation}\label{e:gf-formula}
G_{(0, \infty)}(x, y)
=\int^{x \wedge y}_0 v(z)v(z+|x-y|)dz .
\end{equation}

For any $r>0$, let $G_{(0, r)}$ be the Green function of $Z^{(0,
r)}$,
the subprocess of $Z$ in $(0, r)$. Then
 we have the following result.

\begin{prop}\label{p:upbdongfofkpinfiniteinterval2}
 There exists $c=c(\alpha)>0$ such that for every $r\in (0,
\infty)$,
$$
\int^r_0G_{(0, r)}(x, y)dy\le c(r \wedge r^{\alpha/2})\left((x
\wedge x^{\alpha/2})\wedge((r-x) \wedge (r-x)^{\alpha/2})\right),
\quad x\in (0, r).
$$
\end{prop}

\pf For any $x\in (0, r)$, by \eqref{e:gf-formula}, we have
\begin{eqnarray*}
 \int^r_0G_{(0, r)}(x, y)dy
&\le& \int^r_0G_{(0, \infty)}(x, y)dy\\
&=& \int^x_0\int^x_{x-y}v(z)v(y+z-x)dzdy+
\int^r_x\int^x_0v(z)v(y+z-x)dzdy\\
&=& \int^x_0v(z)\int^x_{x-z}v(y+z-x)dydz
+\int^x_0v(z)\int^r_xv(y+z-x)dydz \\
&\le&   2\,V((0, r))\,V((0, x)).
\end{eqnarray*}
Thus, by \eqref{e:estimates-for-v}
$$
\int^r_0G_{(0, r)}(x, y)dy\le c(r \wedge r^{\alpha/2})(x \wedge
x^{\alpha/2}), \qquad x\in (0, r).
$$
Now the proposition follows by the symmetry.
\qed

Now we return to the process $X^1$ in $\R^d$.
Recall that $C^{\infty}_c(\bR^d)$ is contained in the domain of the
$L_2$-generator $\Delta + \Delta^{\alpha/2}$ of $X^1$ and
$$
(\Delta + \Delta^{\alpha/2}) \phi(x)= \Delta\phi(x)+ \int_{\bR^d}
(\phi(x+y)-\phi(x)-(\nabla \phi(x)\cdot y)1_{B(0, \eps)}(y))
j^1(|y|) dy, \quad\forall  \phi \in C^{\infty}_c(\bR^d)
$$
(see   \cite[Section 4.1]{Sk}). Using the argument in \cite[pp.
152]{KS1}, one can easily see that the last formula on \cite[pp.
152]{KS1} is valid for $X^1$  for  all  $d\ge 1$. Thus  we  have
the following generalization of Dynkin's formula: for every $\phi$
in $C^{\infty}_c(\bR^d)$ and $x\in U$,
\begin{equation}\label{har_gen}
\E_x\left[\phi\left(X^1_{\tau^1_U}\right)\right]-\phi(x) = \int_U G^1_U(x,y)
(\Delta + \Delta^{\alpha/2})  \phi(y)dy=
\E_x \int_0^{\tau^1_U} (\Delta + \Delta^{\alpha/2})  \phi(X^1_s) ds .
\end{equation}

The following
estimates on harmonic measures will play a crucial role in Section
\ref{sec:upper-half-space}.

\begin{thm}\label{T: harmonic}
 For any $R>0$, there exists a constant $c=c(\alpha, R)>0$
such that for every $r \ge R$ and open set  $U \subset B(0,r)$,
$$
 \P_x\left(X^1_{\tau^1_U} \in B(0, r)^c\right)\, \le\, c\,
r^{-\alpha}\int_U  G^1_U(x,y)dy, \quad \text{for every }  x\in U
\cap B(0,r/2).
$$
\end{thm}

\pf Without loss of generality, we assume that $R \in (0, 1)$.
Take a sequence of radial functions $\phi_k$ in $C^{\infty}_c(\bR^d)$
such that $0\le \phi_k\le 1$,
$$
\phi_k(y)=\left\{
\begin{array}{lll}
0, &\hbox{if }  |y|<1/2\\
1, &\hbox{if } 1\le |y|\le k+1\\
0, &\hbox{if } |y|>k+2,
\end{array}
\right.
$$
and that $\sum_{i, j}|\frac{\partial^2}{\partial y_i\partial y_j}
\phi_k|$ is uniformly bounded. Define $\phi_{k,
r}(y)=\phi_k(\frac{y}{r})$. Then we have $0\le \phi_{k, r}\le 1$,
$$
\phi_{k, r}(y)=\left\{ \begin{array}{lll}
0, &\hbox{if } |y|<r/2\\
1, &\hbox{if } r\le |y|\le r(k+1)\\
0, &\hbox{if } |y|>r(k+2),
\end{array} \right.
\qquad \text{and} \qquad
\sup_{y\in \bR^d}
\sum_{i, j}\left|\frac{\partial^2}{\partial y_i\partial y_j}
\phi_{k, r}(y)\right| \,<\, c_{1}\, r^{-2}.
$$
Using this inequality,  we have for $r \ge R$
\begin{align}
&\sup_{k \ge 1}  \sup_{z\in \bR^d} \left|(\Delta +
\Delta^{\alpha/2})\phi_{k,r}(z)\right|\,\le \,
\sup_{k \ge 1}  \sup_{z\in \bR^d} |\Delta
\phi_{k,r}(z)|
+\sup_{k \ge 1}  \sup_{z\in \bR^d} |
\Delta^{\alpha/2}\phi_{k,r}(z)|
 \nonumber\\
&\le  c_{1}\, r^{-2} + \sup_{k
\ge 1}  \sup_{z\in \bR^d} \left|\int_{\bR^d}
(\phi_{k,r}(z+y)-\phi_{k,r}(z)-(\nabla \phi_{k,r}(z)
\cdot y)1_{B(0, r)}(y)) j^1(|y|)dy \right| \nonumber\\
&\le c_{1}\, r^{-2} + c_{2} \sup_{k \ge 1}  \sup_{z\in \bR^d}
\left(\int_{\{|y|\le r\}} \left|
\frac{\phi_{k,r}(z+y)-\phi_{k,r}(z)-(\nabla \phi_{k,r}(z)\cdot y)}
{|y|^{d+\alpha}}\right| dy+\int_{\{r<|y|\}} |y|^{-d-\alpha} dy
\right)\nonumber\\
&\le c_{1}\, r^{-2} +c_{3}\left(\frac{1}{r^2}\int_{\{|y|\le r \}}
\frac{|y|^2}{|y|^{d+\alpha}}dy+ \int_{\{r<|y|\}} |y|^{-d-\alpha} dy
\right)
\,\le \, c_{1}\, r^{-2} + c_{4}r^{-\alpha}.\label{e2.1}
\end{align}
When $U \subset B(0,r)$ for some $r \ge R$, we get, by combining
(\ref{har_gen}) and  (\ref{e2.1}), that for any  $ x\in U \cap
B(0,r/2)$,
$$
\P_x\left(X^1_{\tau^1_U} \in B(0, r)^c\right)    \le \lim_{k\to
\infty} \E_x\left[\phi_{k, r} \left(X^1_{\tau^1_U}\right)\right]
\le c_{5} r^{-\alpha}\int_U  G^1_U(x,y)dy.
$$
\qed

In the remainder of this section we will establish a result (Lemma
\ref{l:lowkey2}) that will be crucial for our heat kernel estimates
in Section \ref{sec:lower-half-space}.

Let
\begin{equation}\label{e:1.1n}
\wh \Delta^{\alpha /2} u(x)\, :=\, \lim_{\eps \downarrow 0}\int_{\{y\in
\bR^d: \, |y-x|>\eps\}} (u(y)-u(x)) \frac{\sA (d, \alpha)
}{|x-y|^{d+\alpha}}\, dy.
\end{equation}
Recall that $\wh \Delta^{\alpha /2}=\Delta^{\alpha /2}$ on $C^\infty_c(\bR^d)$.
 For $x\in \bR^d$ and $p>0$, set
$w_p(x): = (x_d^+)^{p}$. For $0<p<\alpha<2$, let
\begin{equation}\label{e:Lambda}
\Lambda=\Lambda(\alpha,p)
= \frac{p\mathcal{A}(d,-\alpha)}{\alpha }
 \int_0^1\frac{t^{\alpha-p-1}-t^{p-1}}{(1-t)^{\alpha}}dt
\int_{|y|=1,y_d\geq 0}y_d^{\alpha}\  m(dy),
\end{equation}
with the convention that $ m(dy) $ is the Dirac measure when $d=1$.
Then it follows from \cite[Lemma 6.1]{G} that
\begin{align}\label{lag}
\wh \Delta^{\alpha /2} w_p(x)=&\Lambda(d,\alpha,p)w_{p-\alpha}(x),\ \ \ \
x\in \bH.
\end{align}
In particular, on $\bH$ we have
\begin{align} \label{PPPPP}
 \wh \Delta^{\alpha /2} w_p<0,\ \ 0<p<\alpha/2;\ \ \wh
\Delta^{\alpha /2} w_p=0,\ \ p=\alpha/2;\ \ \wh \Delta^{\alpha /2}
w_p>0,\ \ \alpha/2<p<\alpha.
\end{align}

\begin{lemma}\label{l:lowkey1new}
Suppose $  0 < p  \le \frac{\alpha}2$ and $R>8$. Let $Q(a,b):=\{y
\in \bH: |\wt y| <a,  0 <y_d <b\}$ and
$$
h_p(y):=w_p(y) {\bf 1}_{Q(R, R)}(y),  \quad y \in \bH.
$$
There exist constants $c_1, c_2>0$ such that for every $R>8$
and $x \in  Q(2R/3, 2R/3)$,
\begin{equation}\label{e:kl1}
-c_1 (x_d)^{p-\alpha} \leq
\wh \Delta^{\alpha /2} h_p(x) \le -\Lambda (x_d)^{p-\alpha}
\qquad \hbox{when }  0  < p <\frac{\alpha}2
\end{equation}
and
\begin{equation}\label{e:kl2}
-c_1 R^{-\alpha/2} \leq
 \wh \Delta^{\alpha /2} h_{\alpha/2}(x)
\leq  -c_2  R^{-\alpha/2} \qquad \hbox{when } p=\frac{\alpha}2,
\end{equation}
where $\Lambda= \Lambda( \alpha,p)>0$ is the constant defined in
\eqref{e:Lambda}.
\end{lemma}

\pf Since $h_p (y)=w_p(y)$ for  $y \in Q(R, R)$, by \eqref{PPPPP},
we have for any  $x\in Q(2R/3, 2R/3)$,
\begin{align*}
\wh \Delta^{\alpha /2} h_p(x)
&=\wh \Delta^{\alpha /2} (h_p-w_p)(x) +\wh \Delta^{\alpha /2} w_p(x)\\
&= - \int_{Q(R, R)^c}  (y_d^+)^p
\frac{\sA (d,
 -\alpha) }{|x-y|^{d+\alpha}}\, dy+\wh \Delta^{\alpha /2} w_p(x).
\end{align*}
Observe that for $x\in Q(2R/3, 2R/3)$ and $y\in Q(R, R)^c$, $
|y-x|\ge  |y|/3 $.
Thus for $x\in Q(2R/3, 2R/3)$, by the change of variable
$z=R^{-1}y$,
\begin{eqnarray*}
\int_{Q(R, R)^c} \frac{(y_d^+)^p}{|x-y|^{d+\alpha}} dy \le
c_1 \int_{\{y\in\bR^d: \, |y|>R\}} \frac{1}{|y|^{d+\alpha-p}} dy
\le c_2 R^{p-\alpha}\int_{\{z\in \bR^d: \, |z|>1\}}
\frac{1}{|z|^{d+\alpha-p}}\, dz \leq c_3 R^{p-\alpha}.
\end{eqnarray*}
The conclusion of the lemma now follows from the above two displays
and \eqref{lag} and \eqref{PPPPP}.
\qed

\begin{lemma}\label{l:lowkey2}
 There exist
 $c=c(\alpha)>0$ and $R_1=R_1(\alpha)>2$
  such that for every $R>8 R_1$ and
$x\in Q(R/4, R/2) \setminus Q(R/4, 2 R_1)$, we have
$$
\P_x\left(X^1_{\tau^1_{V_R}} \in Q(R, R) \setminus Q(R, R/2)\right)
\ge c \frac{\delta_{\bH} (x)^{\alpha/2}}{R^{\alpha/2}},
$$
where $V_R:=Q(R/2, R/2) \setminus Q(R/2, R_1)$.
\end{lemma}

\pf
Put $p:=(\alpha/4)\vee (\alpha-1)$ and define
$$
h_p(y):=w_p(y) {\bf 1}_{Q(R, R)}(y) \quad \text{and} \quad
h_{\alpha/2}(y):=w_{\alpha/2}(y) {\bf 1}_{Q(R, R)}(y).
$$
We choose $R_1>2$ large such that
\begin{equation}
\frac{\alpha}{2}(1-\frac{\alpha}{2})(R_1/2)^{\alpha-2}\le |\Lambda|,
\label{e:palpal1}
\end{equation}
where $\Lambda$ is the constant defined in \eqref{e:Lambda}.
Obviously, with the above value of $p$, $\Lambda<0$.
For $R>8R_1$ and
$y \in  Q(2R/3, 2R/3)
\setminus Q(R/3, R_1/2)$
by Lemma \ref{l:lowkey1new}
and using the fact that $0 \vee (\frac{3\alpha}2-2) < p <
\frac{\alpha}2<1$, we obtain
\begin{eqnarray*}
&&(\Delta+\wh \Delta^{\alpha /2}) \Big(h_{\alpha/2}(y)- R_1^{\alpha/2-p}
h_p(y)\Big) \\
&\ge& -\frac{\alpha}{2}(1-\frac{\alpha}{2})(y_d)^{\frac{\alpha}{2}-2}
-c_1 R^{-\alpha/2}
-R_1^{\alpha/2-p}\, p(p-1) (y_d)^{p-2}
 +|\Lambda| R_1^{\alpha/2-p} (y_d)^{p-\alpha}  \\
 &=&  (y_d)^{p-\alpha} \left(|\Lambda| R_1^{\alpha/2-p}+ p(1-p)
R_1^{\alpha/2-p}(y_d)^{\alpha-2} -\frac{\alpha}{2}(1-
\frac{\alpha}{2})  (y_d)^{\frac{3\alpha}2-2-p}\right)
 -c_1 R^{-\alpha/2}\\
& \ge&  (y_d)^{p-\alpha} \left(|\Lambda| R_1^{\alpha/2-p}-\frac{\alpha}{2}(1-
\frac{\alpha}{2})  (R_1/2)^{\frac{3\alpha}2-2-p}\right)
 -c_1 R^{-\alpha/2}.
 \end{eqnarray*}
Now, using \eqref{e:palpal1}, we have, for $y \in  Q(2R/3, 2R/3)
\setminus Q(R/3, R_1/2)$
 \bee \label{e:keyl1} (\Delta+\wh
\Delta^{\alpha /2}) \Big(h_{\alpha/2}(y)-R_1^{\alpha/2-p}
h_p(y)\Big) \ge -c_1 R^{-\alpha/2}.
 \eee
Moreover, for $y \in Q(R, R_1)$,
 \bee\label{e:keyl2} (h_{\alpha/2}-R_1^{\alpha/2-p}
h_p)(y)=y_d^{\alpha/2}(1-(R_1/y_d)^{\alpha/2-p}) \le 0.
 \eee
Let $g$ be a nonnegative smooth radial function with compact support
in $\bR^d$ such that $g(x)=0$ for $|x|>1$ and
$\int_{\bR^d}g(x)dx=1$. For $k\ge 1$, define $g_k(x)=2^{kd}g(2^kx)$.
Define
$$
u_{k}(z):= g_k\ast \Big(h_{\alpha/2}-R_1^{\alpha/2-p}
h_p\Big) (z):=\int_{\bR^d} g_k(y)(h_{\alpha/2}-R_1^{\alpha/2-p}
h_p)(z-y)dy \in C^\infty_c(\bR^d).
$$
Let $Q_{R, k}:=\{z\in \bH: {\rm dist}(z, \, Q(R, R))<2^{-k}\}$ and
$A_k=\{x\in \bH: x_d\in ( R_1-2^{-k}, R_1] \}$.
 Note that $u_k=0$ on $Q_{R, k}^c$ and by \eqref{e:keyl2},
for $k$ sufficiently large so that $2^{-k}<R_1/3$,
 \bee \label{e:keyl3}
u_{k}(z)\le 0 \qquad \hbox{for } z_d \leq R_1-2^{-k} ,
 \eee
and for $z\in V_R$, by \eqref{e:keyl1},
 \bee\label{e:keyl4}
(\Delta+\Delta^{\alpha/2})u_{k}(z)=(\Delta+\wh
\Delta^{\alpha/2})u_{k}(z) =g_k\ast (\Delta+\wh
\Delta^{\alpha/2})(h_{\alpha/2}-R_1^{\alpha/2-p} h_p)(z)\ge -c_1
R^{-\alpha/2}.
 \eee
Therefore, using \eqref{har_gen} and
\eqref{e:keyl2}--\eqref{e:keyl4},
 we have that, for any $x\in V_R$,
\begin{align*}
& u_{k}(x)
\,=
\, -\E_x\left[\int_0^{\tau^1_{V_R}}
 (\Delta+\Delta^{\alpha/2})
u_{k}(X^1_t)dt\right]+ \E_x\left[u_{k}
\Big(X^1_{\tau^1_{V_R}}\Big)\right]\\
&\le\, c_1 R^{-\alpha/2} \E_x[\tau^1_{V_R}]+
\E_x\left[u_{k}\Big(X^1_{\tau^1_{V_R}}\Big):X^1_{\tau^1_{V_R}}
 \in Q_{R, k} \setminus Q(R, R_1)\right]
  + \E_x\left[u_{k}\Big(X^1_{\tau^1_{V_R}}\Big):
X^1_{\tau^1_{V_R}} \in A_k \right] \\
&\le\, c_1 R^{-\alpha/2} \E_x[\tau^1_{V_R}]+ \sup_{z\in
A_{k}}|u_{k}(z)|\
\P_x\left( X^1_{\tau^1_{V_R}} \in A_{k}\right)\\
&\quad +\left(\sup_{z \in Q_{R, k} \setminus Q(R,
R_1) } u_{k}(z) \right)\P_x\left(X^1_{\tau^1_{V_R}} \in Q_{R, k}\setminus Q(R,
R_1)\right)\\
&\le\, c_1 R^{-\alpha/2} \E_x[\tau^1_{V_R}]+ \sup_{z\in
A_{k}}|u_{k}(z)| +\left(\sup_{z \in Q_{R, k} } h_{\alpha/2}(z)
\right)\P_x\left(X^1_{\tau^1_{V_R}} \in Q_{R, k} \setminus Q(R,
R_1)\right)\\
&\le\,c_1 R^{-\alpha/2} \E_x[\tau^1_{V_R}]+ \sup_{z\in
A_{k}}|u_{k}(z)|+ R^{\alpha/2}\P_x\left(X^1_{\tau^1_{V_R}}
\in Q_{R, k} \setminus Q(R, R_1)\right).
\end{align*}
Since $h_{\alpha/2}(z)-R_1^{\alpha/2-p} h_p (z)=0$ on $z_d=R_1$,
$ \lim_{k\to \infty}\sup_{z\in A_{k}}|u_{k}(z)|=0$.
Observe that
  $  Q_k(R, R)\setminus Q(R, R_1))$ decreases to $\overline {Q(R,
R)}\setminus Q(R, R_1)$ as $k\to \infty$. We have
\begin{eqnarray*}
\lim_{k\to \infty}\P_x\left(X^1_{\tau^1_{V_R}} \in
Q_{R, k}\setminus Q_1(R, R_1) \right)
&=&\P_x\left(X^1_{\tau^1_{V_R}}\in \overline {Q(R, R)}\setminus Q(R,
R_1)\right)\\
&=& \P_x\left(X^1_{\tau^1_{V_R}}\in Q(R, R) \setminus Q(R, R_1)\right),
\end{eqnarray*}
where the last equality is due to an application of L\'evy system and
the fact that $\partial Q(R, R)$ has zero Lebesgue measure.
Therefore   for $x\in Q(R/2, R/2) \setminus Q(R/2,
2R_1)$, since $x_d \ge 2R_1$,
\begin{eqnarray*}
 (1-2^{p-\alpha/2}) (x_d)^{\alpha/2} &\le&
(x_d)^{\alpha/2}(1-(R_1/x_d)^{\alpha/2-p}) =\lim_{k\to\infty}u_{k}(x)\\
 &\le&
c_1R^{-\alpha/2} \E_x[\tau^1_{V_R}]+
R^{\alpha/2}\P_x\left(X^1_{\tau^1_{V_R}} \in Q(R, R) \setminus Q(R,
R_1)\right),
\end{eqnarray*}
which implies
 \bee
 (x_d)^{\alpha/2} \le  c_1\frac{R^{-\alpha/2} } {1-2^{p-\alpha/2}}
\E_x[\tau^1_{V_R}]+ \frac{R^{\alpha/2}}{1-2^{p-\alpha/2}}
\P_x\left(X^1_{\tau^1_{V_R}} \in Q(R, R) \setminus Q(R, R_1)\right).
\label{e:kubd}
 \eee

Now take a non-negative function $\phi$ in $C^{\infty}_c(\bR^d)$
such that $0\le \phi\le 1$,
$$
\phi(y)=\left\{
\begin{array}{lll}
0 \quad  &\hbox{if } |\wt y|<1/4 \quad \text{or} \quad |y_d|>2, \\
1  &\hbox{if } 1/2\le |\wt y|\le2 \quad \text{and} \quad |y_d|<1, \\
0 &\hbox{if } |\wt y|>3,
\end{array}
\right.
$$
and that $\sum_{i, j}|\frac{\partial^2}{\partial y_i\partial y_j}
\phi|$ is uniformly bounded. Define $\phi_{
R}(y)=\phi(\frac{y}{R})$. Then we have $0\le \phi_{ R}\le 1$,
\bee\label{e:pdef2}
\phi_{ R}(y)=\left\{ \begin{array}{lll}
0 \quad  &\hbox{if } |\wt y|<R/4 \quad \text{or} \quad |y_d|>2R, \\
1 &\hbox{if }  R/2\le |\wt y|\le 2R \quad \text{and} \quad |y_d|<R, \\
0 &\hbox{if } |\wt y|>3R,
\end{array} \right.
\eee
and
$$
\sup_{y\in \bR^d}
\sum_{i, j}\left|\frac{\partial^2}{\partial y_i\partial y_j}
\phi_{ R}(y)\right| \,<\, c_{2}\, R^{-2}.
$$
Using this inequality,  by the argument leading to \eqref{e2.1}, we
get
$$
\sup_{k \ge 1}  \sup_{z\in \bR^d} \left|(\Delta +
\Delta^{\alpha/2})\phi_{R}(z)\right|\,\le\, c_{2}\, R^{-2} \,+\, c_{3}\,R^{-\alpha}.
$$
Thus, by this and Lemma \ref{l:lowkey1new}, for $R>8R_1$
and $y \in Q(2R/3, 2R/3)$, we obtain
\begin{eqnarray}
(\Delta+\wh \Delta^{\alpha /2}) \Big(h_{\alpha/2}(y) +
\frac{2R^{\alpha/2}}{1-2^{p-\alpha/2}}  \phi_{R}(y)\Big) \le
-\frac{\alpha}{2}(1-\frac{\alpha}{2})(y_d)^{\frac{\alpha}{2}-2} +
c_4 R^{\alpha/2}  R^{-\alpha} \le  c_4 R^{-\alpha/2}. \label{e:ketl6}
\end{eqnarray}

For any $k\ge 1$, define
$$
v_{k}(z):=g_k\ast \Big(h_{\alpha/2}
+ \frac{2R^{\alpha/2}}{1-2^{p-\alpha/2}}  \phi_{R}\Big)(z)
\in C^\infty_c(\bR^d).
$$
Put $\Omega_{R}:= Q(R, R/2) \setminus (Q(R,
R_1)    \cup Q(R/2, R/2))$.
By
 \eqref{e:ketl6}, we have
$(\Delta+\Delta^{\alpha/2}) v_{k}(y)  \le  c_4 R^{-\alpha/2}$
for all $
y \in V_R.$
Thus, using
this and \eqref{har_gen},
we have that
for any $k\ge 1$ and $x\in Q(R/4, R/2) \setminus
Q(R/4, 2 R_1)$
\begin{eqnarray*}
 v_{k} (x)
&=& -\E_x\left[\int_0^{\tau^1_{V_R}} (\Delta+\Delta^{\alpha/2})
v_{k}
(X^1_t)dt\right]+ \E_x \left[v_{k}
\Big(X^1_{\tau^1_{V_R}}\Big)\right]\\
&\ge& -c_4 R^{-\alpha/2} \E_x[\tau^1_{V_R}]+
\E_x\left[v_{k}\Big(X^1_{\tau^1_{V_R}}\Big): X^1_{\tau^1_{V_R}} \in
\Omega_{R} \right ].
\end{eqnarray*}
Letting $k\to\infty$ and using
\eqref{e:pdef2}, we get that for any $x\in Q(R/4, R/2) \setminus
Q(R/4, 2 R_1)$ (where $\phi_{R}(x)=0$),\
\begin{eqnarray}
&&(x_d)^{\alpha/2}= \Big(h_{\alpha/2} + \frac{2R^{\alpha/2}}
{1-2^{p-\alpha/2}}  \phi_{R}\Big)(x)=\lim_{k\to \infty}v_{k} (x)\nn\\
&&\ge -c_4 R^{-\alpha/2} \E_x[\tau^1_{V_R}]+ \E_x\left[\Big(h_{\alpha/2}
+ \frac{2R^{\alpha/2}}{1-2^{p-\alpha/2}}  \phi_{R}\Big)
\Big(X^1_{\tau^1_{V_R}}\Big): X^1_{\tau^1_{V_R}}
\in \Omega_{R} \right ]\nn\\
&&
\ge -c_4 R^{-\alpha/2} \E_x[\tau^1_{V_R}] +
\frac{2R^{\alpha/2}}{1-2^{p-\alpha/2}}\P_x\left(X^1_{\tau^1_{V_R}}
\in \Omega_{R} \right)
.\label{e:kubd2}
\end{eqnarray}
Combining \eqref{e:kubd} and \eqref{e:kubd2}, we get
\begin{eqnarray*}
(x_d)^{\alpha/2} &\le&  \frac{c_1R^{-\alpha/2}} {1-2^{p-\alpha/2}}
\E_x[\tau^1_{V_R}]+ \frac{R^{\alpha/2}}{1-2^{p-\alpha/2}}
\P_x\left(X^1_{\tau^1_{V_R}} \in Q(R, R) \setminus Q(R, R_1)\right)\\
&=&  \frac{c_1R^{-\alpha/2}} {1-2^{p-\alpha/2}}  \E_x[\tau^1_{V_R}]
+ \frac{R^{\alpha/2}}{1-2^{p-\alpha/2}}\P_x\left(
X^1_{\tau^1_{V_R}}\in Q(R, R) \setminus Q(R,
R/2)\right)
\\&&+ \frac{R^{\alpha/2}}{1-2^{p-\alpha/2}}\P_x\Big(
X^1_{\tau^1_{V_R}} \in \Omega_{R}\Big)\\
&\le &  \frac{c_1R^{-\alpha/2}} {1-2^{p-\alpha/2}}
\E_x[\tau^1_{V_R}] + \frac{R^{\alpha/2}}{1-2^{p-\alpha/2}}
\P_x\left(X^1_{\tau^1_{V_R}} \in Q(R, R) \setminus Q(R,
R/2)\right)
\\&&+ \frac12\left(c_4 R^{-\alpha/2}
\E_x[\tau^1_{V_R}] + (x_d)^{\alpha/2}\right).
\end{eqnarray*}
Therefore, we conclude that
\begin{equation}\label{e:2.21}
(x_d)^{\alpha/2} \le    \left(\frac{2 c_1} {1-2^{p-\alpha/2}}+ c_4
\right)  R^{-\alpha/2} \E_x[\tau^1_{V_R}] +
\frac{2R^{\alpha/2}}{1-2^{p-\alpha/2}}\P_x\left(X^1_{\tau^1_{V_R}}
\in Q(R, R) \setminus Q(R,
R/2)\right).
\end{equation}
On the other hand, by the L\'evy system of $X^1$,
\begin{eqnarray*}
&&\P_x\left(X^1_{\tau^1_{V_R}} \in Q(R, R) \setminus Q(R, R/2)\right)
\geq  \P_x\left(X^1_{\tau^1_{V_R}} \in Q(R, R) \setminus Q(R, 3R/4)\right) \\
&&= \E_x \left[ \int_0^{\tau^1_{V_R}} \left(\int_{Q(R, R) \setminus
Q(R, 3R/4)} J^1(X_s^1, z) dz  \right)ds \right] \geq c_5 R^{-\alpha}
\, \E_x[\tau^1_{V_R}].
\end{eqnarray*}
This together with \eqref{e:2.21} establishes the lemma.
\qed

\section{Upper bound heat kernel estimates on half-space}
\label{sec:upper-half-space}

In this section we will establish the desired large time upper bound
for $p^1_{\bH}(t, x, y)$.

\begin{lemma}\label{l:ist}
 For any $t_0>0$ and $R>0$, there exists $c=c(\alpha, t_0, R)>1$ such
that for $t \ge t_0$ and $x \in \bH$ with $\delta_{\bH}(x) =x_d  \ge
R$, we have
$$
 \P_{x}(\tau^1_{\bH}>t) \le c \left(\frac{\delta_{\bH}(x)^{\alpha/2}}{\sqrt
t} \wedge 1 \right).
$$
\end{lemma}
\pf Clearly, we can assume $R \le t_0^{1/\alpha}$ and  we only
need to show the theorem for $R \le \delta_H(x) < t^{1/\alpha}$.
Let $u(x)=(x_d^+)^{\alpha/2}+1$ and
$
U(r):=\{x \in \bH; x_d <r\}$. By
\eqref{PPPPP}, for every $x \in \bH$ with $\delta_{\bH}(x) \ge R$,
$$
(\Delta+\wh \Delta^{\alpha /2})u(x) =-\frac{\alpha}{2}
(1-\frac{\alpha}{2})(x_d)^{\alpha/2-2} <0.
$$
Using the same approximation argument as in the proof of Lemma \ref{l:lowkey2} with
$u_{k}(z):=(g_k\ast u)(z)$
where $g_k$ is the function defined in the proof of Lemma \ref{l:lowkey2}
and letting $k \to \infty$, we see that
for $x \in \bH$ with $r >\delta_{\bH}(x) =x_d >R$,
$$
(1+R^{-\alpha/2}) x_d^{\alpha/2} \ge x_d^{\alpha/2}+1=u(x)\ge
\E_x\left[u\left(X^1_{\tau^1_{U(r)}}\right)\right] \ge r^{\alpha/2} \P_x \left(
X^1_{\tau^1_{U(r)}} \in \bH\setminus U(r)\right).
$$
Applying this and Proposition \ref{p:upbdongfofkpinfiniteinterval2},
we get that for $R <\delta_{\bH}(x) < t^{1/\alpha}$.
\begin{eqnarray*}
\P_{x}\left(\tau^1_{\bH}>t\right) &\le& \P_x\left(\tau^1_{U(t^{1/\alpha})} >t\right) +
\P_x \left( X^1_{\tau^1_{U(t^{1/\alpha})}} \in \bH\setminus U(t^{1/\alpha})\right)\\
&\le& \frac1{t} \E_x\left[\tau^1_{U(t^{1/\alpha})}\right] +
(1+R^{-\alpha/2})
\frac{\delta_{\bH}(x)^{\alpha/2}}{\sqrt t} \\
&\le&
c_1 \frac1{t} (t^{1/\alpha} \wedge
t^{1/2})(\delta_{\bH}(x)^{\alpha/2} \wedge \delta_{\bH}(x)) +
(1+R^{-\alpha/2})\frac{\delta_{\bH}(x)^{\alpha/2}}{\sqrt t} \,\le\,
c_2\frac{\delta_{\bH}(x)^{\alpha/2}}{\sqrt t}.
\end{eqnarray*}
\qed

\begin{lemma}\label{l:u_near}
 For every $t_0$ and $R >0,$ there exists $c=c( \alpha, t_0,
R)>1$ such that for every  $(t,x,y) \in [t_0, \infty) \times \bH
\times \bH$ with  $\delta_{\bH}(x)  \ge R$,
$$
 p^1_{\bH}(t,x,y) \le c t^{-d/\alpha
}\left(\frac{\delta_{\bH}(x)^{\alpha/2}}{\sqrt t} \wedge 1 \right) .
$$
\end{lemma}

\pf By the semigroup property and symmetry,
\begin{eqnarray*}
 p^1_{\bH}(t,x,y) &=&\int_{\bH} \int_{\bH} p^1_{\bH}(t/3,x,z)
p^1_{\bH}(t/3,z,w) p^1_{\bH}(t/3,w,y) dzdw\\
&\le &
 \left(\sup_{z,w \in  \R^d} p^1(t/3,z,w) \right)\P_{x}(\tau^1_{\bH}>t/3)
\P_{y}(\tau^1_{\bH}>t/3).
\end{eqnarray*}
Now the lemma follows from Theorem \ref{T:1.1} and Lemma \ref{l:ist}.
\qed

The next lemma and its proof are given in \cite{CKS5} (also see
\cite[Lemma 2]{BGR} and \cite[Lemma 2.2]{CKS}).

\begin{lemma}\label{l:gen1}
Suppose that $U_1,U_3, E$ are open subsets of $\bR^d$ with $U_1,
U_3\subset E$ and ${\rm dist}(U_1,U_3)>0$. Let $U_2 :=E\setminus
(U_1\cup U_3)$.  If $x\in U_1$ and $y \in U_3$, then for all $t >0$,
\begin{equation}\label{eq:ub}
p^1_{E}(t, x, y) \le \P_x\left(X^1_{\tau^1_{U_1}}\in U_2\right)
\left(\sup_{s<t,\, z\in U_2} p^1_E(s, z, y)\right)+ \E_x
\left[\tau^1_{U_1}\right] \left(\sup_{u\in U_1,\, z\in U_3}J^1(u,z)\right).
\end{equation}
\end{lemma}

\begin{lemma}\label{l:u_off1}
 Suppose that $t_0, R>0$. There exists
 $c=c(\alpha, t_0, R)>0$ such that  for every $(t,x,y) \in [t_0, \infty)
\times \bH \times \bH$ with $\delta_{\bH}(x) \ge R$,
$$
 p^1_{\bH}(t,x,y) \le c \left(\frac{\delta_{\bH}(x)^{\alpha/2}
}{\sqrt t} \wedge 1 \right) \left( t^{-d/\alpha}\, \wedge
\frac{t}{|x-y|^{d+\alpha}}  \right) .
$$
\end{lemma}

\pf By Theorem \ref{T:1.1},
Proposition \ref{C:1.2} and Lemma
\ref{l:u_near},  without loss of generality we can assume
$R=t_0^{1/\alpha}$ and it is enough to prove the lemma for
$t_0^{1/\alpha} \le \delta_{\bH}(x) \le (16)^{-1}t^{1/\alpha}$ and
$|x-y| \ge t^{1/\alpha}$. Let $x_0 =(\wt x, 0)$, $U_1:=B( x_0,
8^{-1}t^{1/\alpha}) \cap \bH$, $U_3:= \{z\in \bH: |z-x|>|x-y|/2\}$
and $U_2:=\bH\setminus (U_1\cup U_3)$.

Let $X^1=(X^{1,1}, \dots, X^{1,d})$ and, for any open interval
$(\beta, \gamma)$ in $\R$, let $\wh \tau_{(\beta, \gamma)}:=
\inf\{t>0: X^{1,d} \notin (\beta, \gamma)\}$. Note that, by
Proposition \ref{p:upbdongfofkpinfiniteinterval2} and the assumption
that $16^{-1} t^{1/\alpha}\ge\delta_{\bH}(x)=x_d \ge t_0^{1/\alpha}$, we
have
\begin{equation}\label{e:upae1}
\E_{x}[\tau^1_{U_1}]\,\le\, \E_{x_d}[\wh \tau_{(0, t^{1/\alpha})}]
\,\le\, c_1\, \sqrt t\, x_d^{\alpha/2}\,=  \,c_1\, \sqrt t \,
\delta_{\bH}(x)^{\alpha/2}.
\end{equation}

Since
$$
|z-x| > \frac{|x-y|}2  \ge \frac{1}2 t^{1/\alpha} \quad \text{for
}z\in U_3,
$$
$U_1 \cap U_3 = \emptyset$ and, if $u\in U_1$ and $z\in U_3$, then
\begin{eqnarray}\label{e:n001}
|u-z| \ge |z-x|-|x_0-x|-|x_0-u| \ge |z-x|- 4^{-1} t^{1/\alpha} \ge
\frac{1}{2}|z-x| \ge \frac{1}{4}|x-y|.
\end{eqnarray}
Thus,
\begin{eqnarray}
\sup_{u\in U_1,\, z\in U_3}J^1(u,z) \le \sup_{(u,z):|u-z| \ge
\frac{1}{4}|x-y|}J^1(u,z)  \, \le \, c_3 |x-y|^{-d -\alpha} .
\label{e:n01}
\end{eqnarray}
If $z \in U_2$,
\begin{equation}\label{e:one2}
\frac32 |x-y| \ge |x-y| +|x-z| \ge  |z-y| \ge |x-y| -|x-z| \ge
\frac{|x-y|}2 \ge 2^{-1} t^{1/\alpha}.
\end{equation}
By Theorem \ref{T:1.1} and \eqref{e:one2},
\begin{eqnarray}
\sup_{s\le t,\, z\in U_2} p^1(s, z, y)
&\le& c_4 \sup_{s\le t,\, |z-y| \ge |x-y|/2} \big(s J^1(z,y)\big) +
c_4 \sup_{s\le t,\,  s^{1/2}\ge|z-y| \ge |x-y|/2, } s^{-d/2}\nn\\
&& + c_4 \sup_{s\le t,\, s^{1/2}\le |z-y|,
\atop
\ 1\ge|z-y| \ge
|x-y|/2, } s^{-d/2}e^{-c_5|z-y|^2/s} \nn\\
&\le& c_6 t |x-y|^{-d -\alpha} + 2^{d+\alpha}c_4 \left(\sup_{s\le
t } \frac{s^{\alpha/2}}{|x-y|^{d+\alpha}} \right)\nn\\
&& + c_4    \left(\sup_{a \ge 1} a^{-d/2} e^{-c_5a}\right)
\sup_{1\ge|z-y| \ge |x-y|/2}
|z-y|^{-d}\nn\\
&\le& c_7 t |x-y|^{-d -\alpha} + c_8
\sup_{1\ge|z-y| \ge |x-y|/2}
\frac{|z-y|^{\alpha}}{|x-y|^{d+\alpha}} \le c_9 t |x-y|^{-d
-\alpha}.\label{e:n02}
\end{eqnarray}
Applying Lemma \ref{l:gen1}, \eqref{e:upae1}, \eqref{e:n01}  and
\eqref{e:n02}, we obtain,
\begin{eqnarray*}
p^1_{\bH}(t, x, y) &\le& c_{10} \E_x[\tau^1_{U_1}]|x-y|^{-d -\alpha} +
c_{11} \P_x\Big (X^1_{\tau^1_{U_1}}\in U_2 \Big)  t |x-y|^{-d -\alpha} \\
&\le& c_{12} \sqrt t \, \delta_{\bH}(x)^{\alpha/2} |x-y|^{-d -\alpha} +
c_{11}\P_x\Big(X^1_{\tau^1_{U_1}}\in U_2\Big)  t |x-y|^{-d -\alpha}.
\end{eqnarray*}
Finally,
applying Theorem \ref{T: harmonic} with $U=U_1$ and
$r=8^{-1}t^{1/\alpha} \ge 2{t_0}^{1/\alpha} $, we have
$$
\P_x \Big(X^1_{\tau^1_{U_1}}\in U_2 \Big) \le \P_x
\Big(X^1_{\tau^1_{U_1}}\in B( x_0, 8^{-1}t^{1/\alpha})^c\Big)
\,\le\, c_{14} \,\frac1t \,\int_{U_1} G^1_{U_1}(x,y)dy \,=\,c_{14}\,
\frac1t\, \E_{x}[\tau^1_{U_1}].
$$
Now applying \eqref{e:upae1}, we have proved
the lemma.
\qed

\begin{lemma}\label{l:up3}
 For every $R>0$ and $t_0>0$, there exists a constant $c=c(R,
\alpha, t_0)$ such that for all $(t, x, y)\in [t_0, \infty) \times
\bH\times \bH$ with $\delta_{\bH}(x) \wedge \delta_{\bH}(y) \ge R$.
\begin{eqnarray*}
 p^1_{\bH}(t, x, y)&\le& c \left(\frac{\delta_{\bH}(x)^{\alpha/2}
}{\sqrt t} \wedge 1 \right) \left(\frac{\delta_{\bH}(y)^{\alpha/2}
}{\sqrt t} \wedge 1 \right)\left( t^{-d/\alpha}\, \wedge
\frac{t}{|x-y|^{d+\alpha}} \right).
\end{eqnarray*}
\end{lemma}

\pf By Lemma \ref{l:u_off1} and Theorem \ref{T:1.1}, we only need to
to prove the theorem for $\delta_{\bH}(x) \vee \delta_{\bH}(y) \le
t^{1/\alpha}$. Denote by $q(t, x, y)$ the transition density of the
$\alpha$-stable process $Y$ in $\R^d$. By Lemma \ref{l:u_off1} and
the lower bound estimate of $q(t,x,y)$,
there is a constant $c_1>0$ so that
$$
p^1_{\bH}(t/2, x, z)\le c_1
\left(\frac{ \delta_{\bH}(x)^{\alpha/2}}{\sqrt{t}} \wedge 1 \right)
q(t/2, x,z) \quad \text{and} \quad
p^1_{\bH}(t/2, z, y)\le c_1
\left(\frac{ \delta_{\bH}(y)^{\alpha/2}}{\sqrt{t}} \wedge 1 \right)
q(t/2, y,z).
$$
Thus, by semigroup property and the upper bound estimate of
$q(t,x,y)$,
\begin{align*}
p^1_{\bH}(t,x,y)
&= \int_{\bH} p^1_{\bH}(t/2,x,z)p^1_{\bH}(t/2,z,y)dz \\
&\le c_2^2 \left(\frac{
\delta_{\bH}(x)^{\alpha/2}}{\sqrt{t}} \wedge 1 \right)
\left(\frac{ \delta_{\bH}(y)^{\alpha/2}}{\sqrt{t}} \wedge
1 \right) \int_{\bH} q(t/2, x,z) q(t/2, y,z)dz\\
&\le c_2^2 \left(\frac{
\delta_{\bH}(x)^{\alpha/2}}{\sqrt{t}} \wedge 1 \right)
\left(\frac{ \delta_{\bH}(y)^{\alpha/2}}{\sqrt{t}} \wedge
1 \right) q(t,x,y)\\
&\le c_3\left(\frac{\delta_{\bH}(x)^{\alpha/2} }{\sqrt t} \wedge 1
\right) \left(\frac{\delta_{\bH}(y)^{\alpha/2} }{\sqrt t} \wedge 1
\right)\left( t^{-d/\alpha}\, \wedge \frac{t}{|x-y|^{d+\alpha}}
\right).
\end{align*}
\qed

\begin{thm}\label{t:ub}
Let $t_0$ be a positive constant.  Then there exists a constant
 $c= c(\alpha,  t_0 )>0$  such that for all $t \in [t_0,
\infty)$ and $x, y \in \bH$,
$$
 p^1_{\bH}(t, x, y)\le c \left(\frac{\delta_{\bH}(x)\wedge
\delta_{\bH}(x)^{\alpha/2}}{\sqrt{t}} \wedge 1 \right)
\left(\frac{\delta_{\bH}(y) \wedge
\delta_{\bH}(y)^{\alpha/2}}{\sqrt{t}} \wedge 1 \right)\left(
t^{-d/\alpha}\, \wedge \frac{t}{|x-y|^{d+\alpha}}  \right).
$$
\end{thm}

\pf Let $x_0$ and $y_0$ be as in \eqref{e:x0y0}. By the semigroup
property and \eqref{E:ratio-half-space-like1},  we have
\begin{align}
p^1_{\bH}(t,x,y)
&= \int_{\bH} \int_{\bH }p^1_{\bH}(t_0,x,z)p^1_{\bH}
(t-2t_0,z,w)p^1_{\bH}(t_0,w,y)dzdw \nn\\
&\asymp \left( 1\wedge \delta_{\bH}(x) \right)\left(
1\wedge \delta_{\bH}(y) \right)
\int_{\bH} \int_{\bH} p^1_{\bH}(t_0,x_0,z)p^1_{\bH}(t-2t_0,z,w)
p_{\bH}(t_0,w,y_0)dzdw
\nn\\
&= \left( 1\wedge \delta_{\bH}(x) \right)\left( 1\wedge \delta_{\bH}(y)
\right)p^1_{\bH}(t,x_0,y_0). \label{e:asdfgh}
\end{align}
By Lemma \ref{l:up3} and the fact $|x_0-y_0|=|x-y|$, we have
\begin{align*}
p^1_{\bH}(t,x_0,y_0) &\le
c_1   \left(\frac{\delta_{\bH}(x_0)^{\alpha/2} }{\sqrt t} \wedge
1 \right)
\left(\frac{\delta_{\bH}(y_0)^{\alpha/2} }{\sqrt t} \wedge
1 \right)\left( t^{-d/\alpha}\, \wedge
\frac{t}{|x-y|^{d+\alpha}}  \right).
\end{align*}
This together with Lemma \ref{l:piw2} (with $a=1$ there)
and \eqref{e:asdfgh} proves the theorem.
\qed

\section{Lower bound heat kernel estimates on half-space}\label{sec:lower-half-space}

In this section we establish the desired sharp large time lower
bound on $p^1_{\bH}(t, x, y)$. We will use some ideas from
\cite{BGR, CKS5}.

\begin{lemma}\label{L:4.2}
 For any positive constant  $t_0$, there exists $c=c(t_0,
\alpha)>0$ such that for any
$t \ge t_0$ and $y \in \R^d$,
$$
\P_y \left( \tau^1_{B(y, 8^{-1} t^{1/\alpha})} > t/3 \right)\,
\ge\, c.
$$
\end{lemma}

\pf By \cite[Proposition 6.2]{CK08}, there exists $\eps=\eps(t_0,
 \alpha)>0$ such that for every $t \ge t_0$,
$$
\inf_{y\in \bR^d } \P_y \left( \tau^1_{B(y, {16}^{-1} t^{1/\alpha}
)}
> \eps t \right) \ge \frac12.
$$
Suppose $\eps <\frac13 $, then by the parabolic Harnack inequality
in \cite{CK08, SV07},
$$
c_1\,p ^1_{B(y,8^{-1} t^{1/\alpha})}(\eps  t ,y,w) \, \le \,
p ^1_{B(y,8^{-1} t^{1/\alpha})}(t/3 ,y,w)\qquad \hbox{for }
w \in B(y, {16}^{-1} t^{1/\alpha} ),
$$
where the constant $c_1=c_1(t_0, \alpha)>0$ is
independent of $y\in \bR^d$. Thus
\begin{eqnarray*}
\P_y \left( \tau^1_{B(y,8^{-1} t^{1/\alpha})} > t/3\right) &=&
\int_{B(y, 8^{-1} t^{1/\alpha} )}
p^1_{B(y, 8^{-1} t^{1/\alpha} )}(t/3 ,y,w) dw\\
&\ge& c_1\int_{B(y, {16}^{-1} t^{1/\alpha} )} p^1_{B(y, 8^{-1} t^{1/\alpha}
)}(\eps t ,y,w) dw  \ge  \frac{c_1}2.
\end{eqnarray*}
\qed

The next result holds for any symmetric discontinuous
Hunt process that possesses a transition density and whose L\'evy system
admitting  jumping density kernel.
 Its proof is the same as that of \cite[Lemma 3.3]{CKS4} and so it is omitted here.

\begin{lemma}\label{l:gen}
Suppose that $U_1,U_2, U$ are open subsets of $\bR^d$ with $U_1,
U_2\subset U$ and ${\rm dist}(U_1,U_2)>0$. If $x\in U_1$ and $y \in
U_2$, then for all $t >0$,
\begin{equation}\label{eq:lb}
p^1_{U}(t, x, y)\,\ge\, t\,  \P_x(\tau^1_{U_1}>t)
\,\P_y(\tau^1_{U_2}>t)\inf_{u\in U_1,\, z\in U_2}J^1(u,z) \,.
 \end{equation}
\end{lemma}

\begin{lemma}\label{lower bound12}
 Suppose that $t_0>0$. There exists $c=c(t_0, \alpha)>0$ such that
for all $t \ge t_0$ and  $u, v\in \bR^d$ with $|u-v|\ge
t^{1/\alpha}/2$,
$$
 p^1_{B(u,t^{1/\alpha})\cup B(v,t^{1/\alpha})}(t/3, u, v)\,\ge\, c \,
t \, |u- v|^{-d-\alpha}.
$$
\end{lemma}

\pf Let $U= B(u,t^{1/\alpha})\cup B(v,t^{1/\alpha})$. With  $U_1=
B(u,t^{1/\alpha}/8)$ and $U_2=B(v,t^{1/\alpha}/8)$, we have by Lemma
\ref{l:gen} that
\begin{eqnarray*}
p^1_{U}(t/3, u, v) \ge \frac{t}3 \P_u(\tau^1_{U_1}>t/3)\left(
\inf_{w\in U_1,\, z\in U_2} j^1(|w-z|)\right) \P_v(\tau^1_{U_2} >
t/3) \,.
\end{eqnarray*}
Moreover, $|w-z| \le |u-v| + |w-u|+ |z-v| \le |u-v| + t^{1/\alpha}/4
\le \frac32 |u-v|$. Thus by Lemma \ref{L:4.2},
\begin{eqnarray*}
p^1_{B(u,t^{1/\alpha})\cup B(v,t^{1/\alpha})}(t/3, u, v)&\ge&
\frac{t}3 \left(\P_0(\tau^1_{B(0,t^{1/\alpha}/8)}>t/3)\right)^2
\left(\inf_{ w\in U_1,\, z\in U_2}j^1(|w-z|)\right)\\
&\ge& c_1 \,t \,|u- v|^{-d-\alpha}.
\end{eqnarray*}
\qed

The next result follows from \cite[Proposition 3.4]{SV07}.

\begin{lemma}\label{T:3.4}
 There exist
$ R_2=R_2(\alpha )>1$ and $c=c(\alpha )>0$
such that for all $t\ge R_2^{\alpha}$,
$$
\inf_{x, y\in B(0, 6t^{1/\alpha})}p^1_{B(0,
12t^{1/\alpha})}(t/3, x, y)
\,\ge\, c\, t^{-d/\alpha}.
$$
\end{lemma}

For the remainder of this section, we define $R_3:=R_1\vee
R_2$,
where
$R_1>0$ is the constant in Lemma \ref{l:lowkey2}.
For any $x\in \R^d$ and $a, b>0$, we define
$$
Q_x(a, b):=\{y\in \bH: |\wt y-\wt x|<a, y_d<b\}.
$$

\begin{lemma}\label{l:keylon}
 There is a positive constant $c=c (\alpha)$ such that for all
$(t, x)\in ((4R_1)^\alpha, \infty) \times \bH$ with $2 R_1 <
\delta_{\bH}(x) < t^{1/\alpha}/2$,
$$
\P_x(\tau^1_{Q_x(2t^{1/\alpha}, 2t^{1/\alpha}) }>t/3) \ge c
\frac{\delta_{\bH}(x)^{\alpha/2}}{\sqrt t}.
$$
\end{lemma}

\pf Without loss of generality we assume that $\wt x=\wt 0$ and
let $Q(a, b):=Q_0(a, b)$.
Let $V(t):= Q( t^{1/\alpha}/2,  t^{1/\alpha}/2) \setminus
Q(t^{1/\alpha}/2, R_1)$.
By Lemma \ref{l:lowkey2},
Lemma \ref{L:4.2} and the strong Markov property,
\begin{eqnarray*}
&& \P_x\left( \tau^1_{Q(2t^{1/\alpha}, 2t^{1/\alpha})} >t/3 \right)\\
&\geq & \P_x\left( \tau^1_{Q(2t^{1/\alpha}, 2t^{1/\alpha})} >t/3, \
X^1_{\tau^1_{V(t)}} \in Q(t^{1/\alpha}, t^{1/\alpha}) \setminus
Q(t^{1/\alpha}, t^{1/\alpha}/2)\right) \\
& = &\E_x\left[ \P_{X^1_{\tau^1_{V(t)
}}}\left(\tau^1_{Q(2t^{1/\alpha}, 2t^{1/\alpha})} >t/3 \right):
X^1_{\tau^1_{V(t) }} \in Q(t^{1/\alpha}, t^{1/\alpha})
\setminus Q(t^{1/\alpha},  t^{1/\alpha}/2)\right]\\
&\geq   & \E_x\left[ \P_{X^1_{\tau^1_{V(t)
}}}
\left(\tau^1_{B(X^1_{\tau^1_{V(t)}}, \, 4^{-1}t^{1/\alpha})}>t/3 \right)
: X^1_{\tau^1_{V(t) }} \in Q(t^{1/\alpha}, t^{1/\alpha})
\setminus Q(t^{1/\alpha},  t^{1/\alpha}/2)\right]\\
&\geq & c_1\P_x\left(X^1_{\tau^1_{V(t) }}\in Q(t^{1/\alpha}, t^{1/\alpha})
\setminus Q(t^{1/\alpha}, 2^{-1}t^{1/\alpha})\right) \ \ge \ c_2
\frac{\delta_{\bH}(x)^{\alpha/2}}{\sqrt t}.
\end{eqnarray*}
This proves the Lemma. \qed

\begin{lemma}\label{l:lowkeyfinal}
 There is a positive constant $c=c(\alpha)$ such that for all
$(t, x,y)\in [(4R_3)^\alpha, \infty) \times \bH\times \bH$ with
$\delta_{\bH}(x) \wedge \delta_{\bH}(y) \ge 2R_3$,
$$
 p^1_{\bH}(t, x, y)\ge c \left(\frac{
\delta_{\bH}(x)^{\alpha/2}}{\sqrt{t}} \wedge 1 \right)
\left(\frac{\delta_{\bH}(y)^{\alpha/2}}{\sqrt{t}} \wedge 1
\right)\left(t^{-d/\alpha} \wedge \frac{ t }{|x-y|^{d+\alpha} }
\right) .
$$
\end{lemma}

\pf Fix $x,y\in \bH$.
Let $x_0=(\wt x, 0)$,  $y_0=(\wt y, 0)$,
$\xi_x:=x+(\wt 0, 32t^{1/\alpha})$ and $\xi_y:=y+(\wt 0,
32t^{1/\alpha})$.
If $2 R_3  \le \delta_{\bH}(x) < t^{1/\alpha}/2$,
by Lemmas \ref{L:4.2}, \ref{l:gen} and  \ref{l:keylon},
\begin{align*}
& \int_{B(\xi_x, 2t^{1/\alpha})} p^1_{{\bH}}(t/3,x,u)du\\
\ge& t\,  \P_x\left(
\tau^1_{Q_x(2 t^{1/\alpha},
2 t^{1/\alpha})}>t/3\right)\left(\inf_{v\in Q_x(2 t^{1/\alpha}, 2 t^{1/\alpha} ) \atop
w\in B(\xi_x,4 t^{1/\alpha})}J^1(v,w) \right)
 \,\int_{B(\xi_x, 2t^{1/\alpha})}\P_u\left(\tau^1_{B(\xi_x, 4t^{1/\alpha})}>
t/3\right)du\\
\ge& c_1 t\,  \P_x\left(\tau^1_{Q_x(2 t^{1/\alpha}, 2 t^{1/\alpha} )
}>t/3\right)
t^{-d/\alpha-1} \,   \P_0\left(\tau^1_{B(0, t^{1/\alpha}/8)}> t/3\right)
\,|B(\xi_x,2 t^{1/\alpha})| \nonumber\\
\ge& c_2 \P_x\left(\tau^1_{Q_x(2 t^{1/\alpha}, 2 t^{1/\alpha} ) }
>t/3\right)\ \ge \ c_3 \frac{\delta_{\bH}(x)^{\alpha/2}}{\sqrt{t}}.
\end{align*}

On the other hand, if $ \delta_{\bH}(x) \ge t^{1/\alpha}/2 \ge
2R_3$, by Lemmas \ref{L:4.2} and \ref{l:gen},
\begin{align*}
& \int_{B(\xi_x, 2t^{1/\alpha})} p^1_{{\bH}}(t/3,x,u)du  \nonumber
\\\ge& t\,  \P_x\left(\tau^1_{B(x,  8^{-1} t^{1/\alpha}) \cap \bH
}>t/3\right)\left(\inf_{v\in B(x_0, 2t^{1/\alpha}) \cap \bH \atop w\in
B(\xi_x,4 t^{1/\alpha})}J^1(v,w) \right) \,\int_{B(\xi_x,
2t^{1/\alpha})}\P_u\left(\tau^1_{B(\xi_x,
4t^{1/\alpha})}>t/3\right)du \nonumber\\
\ge& c_4 t\,  \P_x\left(\tau^1_{B(x,  8^{-1} t^{1/\alpha})}>t/3\right)
t^{-d/\alpha-1}\,  \P_0\left(\tau^1_{B(0, t^{1/\alpha}/8)}>
t/3\right) \,|B(\xi_x,2 t^{1/\alpha})| \nonumber\\
\ge& c_5 \P_x\left(\tau^1_{B(x,  8^{-1} t^{1/\alpha})}>t/3\right)  \ge c_6 .
\end{align*}
Thus
\begin{eqnarray}
\int_{B(\xi_x, 2 t^{1/\alpha})} p^1_{{\bH}}(t/3,x,u)du \ge c_7
\left(1\wedge \frac{\delta_{\bH}(x)^{\alpha/2}}{\sqrt{t}} \right),
\label{e:loww_0}
\end{eqnarray}
and similarly,
\begin{eqnarray}
\int_{B(\xi_y, 2 t^{1/\alpha})} p^1_{{\bH}}(t/3,y,u)du \ge c_7
\left(1\wedge \frac{\delta_{\bH}(y)^{\alpha/2}}{\sqrt{t}} \right).
\label{e:loww_01}
\end{eqnarray}

Now we deal with the cases $|x-y| \ge 5 t^{1/\alpha} $ and $|x-y| <
5t^{1/\alpha}$ separately.

\medskip

\noindent {\it Case 1}: Suppose that $|x-y| \ge  5 t^{1/\alpha} $.
Note that by the semigroup property  and Lemma~\ref{lower bound12},
\begin{align*}
&p^1_{\bH}(t,x,y)\nonumber\\
\geq& \int_{B(\xi_y, 2 t^{1/\alpha})}\int_{B(\xi_x, 2 t^{1/\alpha})}
p^1_{\bH}(t/3,x,u) p^1_{\bH}(t/3,u,v)p^1_{\bH}
(t/3,v,y)dudv \nonumber\\
\geq& \int_{B(\xi_y, 2 t^{1/\alpha})}\int_{B(\xi_x, 2 t^{1/\alpha})}
p^1_{{\bH}}(t/3,x,u)p^1_{B(u, t^{1/\alpha}) \cup
B(v,t^{1/\alpha})}(t/3,u,v)p^1_
{\bH}(t/3,v,y)dudv\nonumber\\
\geq& c_8 t \left(\inf_{(u,v) \in B(\xi_x, 2 t^{1/\alpha}) \times
B(\xi_y,2 t^{1/\alpha})} |u-v|^{-d-\alpha}\right) \int_{B(\xi_y, 2
t^{1/\alpha})}\int_{B(\xi_x, 2 t^{1/\alpha})}
p^1_{{\bH}}(t/3,x,u)p^1_{\bH}(t/3,v,y)dudv.
\end{align*}
It then follows from \eqref{e:loww_0}--\eqref{e:loww_01} that
\begin{align}\label{e:loww1}
p^1_{\bH}(t,x,
y)
\ge
 c_9  t \left(\inf_{(u,v) \in
B(\xi_x, 2 t^{1/\alpha}) \times B(\xi_y,2 t^{1/\alpha})} |u-v|^{-d-\alpha} \right)
\left(\frac{\delta_{\bH}(x)^{\alpha/2}}{\sqrt{t}} \wedge 1 \right)
\left(\frac{\delta_{\bH}(y)^{\alpha/2}}{\sqrt{t}} \wedge 1 \right).
\end{align}
Using  the assumption $|x-y| \ge 5 t^{1/\alpha}$ we get that, for $u\in
B(\xi_x, 2 t^{1/\alpha})$ and $v\in B(\xi_y,2 t^{1/\alpha})$, $|u-v| \le
4 t^{1/\alpha} +|x-y| \le  2 |x-y|$. Hence
\begin{equation}\label{e:loww2}
\inf_{(u,v) \in B(\xi_x, 2 t^{1/\alpha}) \times B(\xi_y,2 t^{1/\alpha})}
|u-v|^{-d-\alpha}\,\ge\, c_{10}
|x-y|^{-d-\alpha}.
\end{equation}
By \eqref{e:loww1} and \eqref{e:loww2}, we conclude that for $|x-y|
\ge 5  t^{1/\alpha}$
\begin{eqnarray*}
p^1_{\bH}(t, x, y)\ge c_{11} \left(\frac{
\delta_{\bH}(x)^{\alpha/2}}{\sqrt{t}} \wedge 1 \right)
\left(\frac{\delta_{\bH}(y)^{\alpha/2}}{\sqrt{t}} \wedge 1 \right)t
|x-y|^{-d-\alpha}.
\end{eqnarray*}

\noindent {\it Case 2}: Suppose $|x-y| < 5 t^{1/\alpha}$. In this
case, for every $(u,v) \in  B(\xi_x, 2 t^{1/\alpha}) \times B(\xi_y,
2 t^{1/\alpha})$, $|u-v| \le 9 t^{1/\alpha}$. Thus, using the fact
that $\delta_{\bH}(\xi_x) \wedge \delta_{\bH}(\xi_y) \ge 32 t^{1/\alpha}$,
there exists $w_0 \in \bH$ such that
\begin{equation}\label{e:dfeggg}
B(\xi_x, 2 t^{1/\alpha}) \cup B(\xi_y, 2 t^{1/\alpha})
\subset B(w_0, 6 t^{1/\alpha})\subset B(w_0,
12 t^{1/\alpha})\subset \bH.
\end{equation}
Now, by the semigroup property and \eqref{e:dfeggg}, we get
\begin{align*}
&p^1_{\bH}(t,x,y)\nonumber\\
\geq& \int_{B(\xi_y, 2 t^{1/\alpha})}\int_{B(\xi_x, 2 t^{1/\alpha})}
p^1_{{\bH}}(t/3,x,u)p^1_{B(w_0, 8t^{1/\alpha})}(t/3,u,v)p^1_
{\bH}(t/3,v,y)dudv\nonumber\\
\geq&  \left(\inf_{u,v \in B(w_0, 6 t^{1/\alpha})} p^1_{B(w_0,
12t^{1/\alpha})}(t/3,u,v)\right) \int_{B(\xi_y, 2
t^{1/\alpha})}\int_{B(\xi_x, 2 t^{1/\alpha})}
p^1_{{\bH}}(t/3,x,u)p^1_{\bH}(t/3,v,y)dudv.
\end{align*}
It then follows from \eqref{e:loww_0}--\eqref{e:loww_01} and
Lemma \ref{T:3.4} that
\begin{align*}
p^1_{\bH}(t, x, y) \ge c_{12} \left(\frac{
\delta_{\bH}(x)^{\alpha/2}}{\sqrt{t}} \wedge 1 \right)
\left(\frac{\delta_{\bH}(y)^{\alpha/2}}{\sqrt{t}} \wedge 1
\right)t^{-d/\alpha}.
\end{align*}

Combining these two cases, we have proved the theorem.
\qed

\begin{thm}\label{t:low}
 There exists a positive constant $c= c(\alpha )$
such that for all $t \in [(4R_3)^\alpha, \infty)$ and $x, y\in \bH$,
$$
 p^1_{\bH}(t, x, y)\ge c \left(1\wedge \frac{\delta_{\bH}(x) \wedge
\delta_{\bH}(x)^{\alpha/2}}{\sqrt{t}} \right)\left(1\wedge
\frac{\delta_{\bH}(y)\wedge \delta_{\bH}(y)^{\alpha/2}}{\sqrt{t}}
\right)\left(t^{-d/\alpha} \wedge \frac{ t }{|x-y|^{d+\alpha} }
\right).
$$
\end{thm}

\pf Let $t_0=(4R_3)^2 > (4R_3)^\alpha$  and let $x_0$ and $y_0$ be
as in \eqref{e:x0y0}. By the semigroup property and
\eqref{E:ratio-half-space-like1} we have
\begin{align}
p^1_{\bH}(t,x,y)
&= \int_{\bH} \int_{\bH} p^1_{\bH}(t_0,x,z)p^1_{\bH}
(t-2t_0,z,w)p^1_{\bH}(t_0,w,y)dzdw \nn\\
&\asymp \left( 1\wedge \delta_{\bH}(x) \right)\left(
1\wedge \delta_{\bH}(y) \right)
\int_{\bH} \int_{\bH} p^1_{\bH}(t_0,x_0,z)p^1_{\bH}(t-2t_0,z,w)
p^1_{\bH}(t_0,w,y_0)dzdw \nn\\
&= \left( 1\wedge \delta_{\bH}(x) \right)\left( 1\wedge \delta_{\bH}(y)
\right)p^1_{\bH}(t,x_0,y_0). \label{e:asfd}
\end{align}
Since, $\delta_{\bH}(x_0) \wedge \delta_{\bH}(y_0)> t_0^{1/2} =4R_3$,
by Lemma \ref{l:lowkeyfinal} and the fact $|x_0-y_0|=|x-y|$,
\begin{align*}
p^1_{\bH}(t,x_0,y_0) &\ge
c_{1} \left(\frac{
\delta_{\bH}(x_0)^{\alpha/2}}{\sqrt{t}} \wedge 1 \right)
\left(\frac{\delta_{\bH}(y_0)^{\alpha/2}}{\sqrt{t}} \wedge 1
\right)\left(t^{-d/\alpha} \wedge \frac{ t }{|x-y|^{d+\alpha} } \right).
\end{align*}
The conclusion of the theorem now follows from the above inequality,
Lemma \ref{l:piw2} and \eqref{e:asfd}.
\qed

\section{Heat kernel estimates on half-space-like
domains}
\label{sec:half-space-like}

In this section, we will establish the main result of this paper.

Combining Theorem \ref{t:main}(i), Theorems \ref{t:ub} and
\ref{t:low}, we get that for every $T>0$, there exist constants
$c_i= c_i(\alpha,  T)\geq 1$, $i=1, 2,$ such that for all $(t, x,
y)\in (0, T]\times \bH\times \bH$,
\begin{eqnarray*}
&&c_1^{-1}\, \left(1\wedge \frac{\delta_{\bH}(x)}{\sqrt{t}}
\right)\left(1\wedge \frac{\delta_{\bH}(y)}{\sqrt{t}} \right) \left(
t^{-d/2} e^{-c_2|x-y|^2/t}+\left( \frac{t}{|x-y|^{d+\alpha}}
\wedge t^{-d/2} \right)\right)\\
&& \le p^1_{\bH}(t, x, y) \leq c_1 \left(1\wedge
\frac{\delta_{\bH}(x)}{\sqrt{t}} \right)\left(1\wedge
\frac{\delta_{\bH}(y)}{\sqrt{t}}  \right) \left( t^{-d/2}
e^{-|x-y|^2/(c_2t)}+\left( \frac{t}{|x-y|^{d+\alpha}} \wedge
t^{-d/2} \right)\right)
\end{eqnarray*}
and for all $t \in [T, \infty)$ and $x, y$ in $H$,
\begin{align*}
&c_1^{-1} \left(1\wedge \frac{\delta_{\bH}(x) \wedge
\delta_{\bH}(x)^{\alpha/2}}{\sqrt{t}} \right)\left(1\wedge
\frac{\delta_{\bH}(y)\wedge \delta_{\bH}(y)^{\alpha/2}}{\sqrt{t}}
\right)\left(t^{-d/\alpha} \wedge \frac{ t }{|x-y|^{d+\alpha} } \right)\\
& \le p^1_{\bH}(t, x, y)\le c_1 \left(1\wedge \frac{\delta_{\bH}(x) \wedge
\delta_{\bH}(x)^{\alpha/2}}{\sqrt{t}} \right)\left(1\wedge
\frac{\delta_{\bH}(y)\wedge \delta_{\bH}(y)^{\alpha/2}}{\sqrt{t}}
\right)\left(t^{-d/\alpha} \wedge \frac{ t }{|x-y|^{d+\alpha} }
\right).
\end{align*}

Now using \eqref{e:scaling}, we established Theorem \ref{t:main_hsl}
for $D=\bH$  in the form of \eqref{e:old_1}--\eqref{e:old_2}.
\begin{thm}\label{t:2sbhalfspace}
 For every $T>0$, there exist $c=c(\alpha,  T)\geq 1$ and $C_3=C_3(\alpha,
 T)\geq 1$ such that for all $a>0$ and $(t, x, y)\in (0,
a^{2\alpha/(\alpha-2)}T]\times \bH\times \bH$,
\begin{eqnarray*}
 &&c^{-1}\, \left(1\wedge \frac{\delta_{\bH}(x)}{\sqrt{t}}
\right)\left(1\wedge \frac{\delta_{\bH}(y)}{\sqrt{t}} \right) \left(
t^{-d/2} e^{-C_3|x-y|^2/t}+\left( \frac{a^\alpha
t}{|x-y|^{d+\alpha}}
\wedge t^{-d/2} \right)\right)\\
 && \le p^a_{\bH}(t, x, y) \leq c \left(1\wedge
\frac{\delta_{\bH}(x)}{\sqrt{t}} \right)\left(1\wedge
\frac{\delta_{\bH}(y)}{\sqrt{t}}  \right) \left( t^{-d/2}
e^{-|x-y|^2/(C_3t)}+\left( \frac{a^\alpha t}{|x-y|^{d+\alpha}}
\wedge t^{-d/2} \right)\right)
\end{eqnarray*}
and for all
$t \in [a^{2\alpha/(\alpha-2)}T, \infty)$
and $x, y$ in $\bH$,
\begin{align*}
 &c^{-1} \left(1\wedge \frac{\delta_{\bH}(x) \wedge
(a^{-1}\delta_{\bH}(x))^{\alpha/2}}{\sqrt{t}} \right)\left(1\wedge
\frac{\delta_{\bH}(y)\wedge
(a^{-1}\delta_{\bH}(y))^{\alpha/2}}{\sqrt{t}} \right)\left((a^\alpha
t)^{-d/\alpha} \wedge \frac{a^\alpha t
}{|x-y|^{d+\alpha} } \right)\\
 & \le p^a_{\bH}(t, x, y)
 \\
 &\le
  c \left(1\wedge \frac{\delta_{\bH}(x) \wedge
(a^{-1}\delta_{\bH}(x))^{\alpha/2}}{\sqrt{t}} \right)\left(1\wedge
\frac{\delta_{\bH}(y)\wedge
(a^{-1}\delta_{\bH}(y))^{\alpha/2}}{\sqrt{t}} \right)\left((a^\alpha
t)^{-d/\alpha} \wedge \frac{a^\alpha t }{|x-y|^{d+\alpha} } \right).
\end{align*}
\end{thm}

Now we are in a position to establish the main result of this paper.

\bigskip

\noindent {\bf Proof of  Theorem \ref{t:main_hsl}}. Recall that that
$D$ is a half-space-like $C^{1,1}$ domain with $C^{1, 1}$
characteristics $(R_0, \Lambda_0)$ and $\bH_b\subset D\subset \bH$
for some $b>0$ such that
 that the path distance in $D$ is comparable to
the Euclidean distance with characteristic $\lambda_0$.
 Then we have the following trivial inequalities
\begin{equation}\label{e:bdbyhs}
p^a_{\bH_b}(t, x, y)\le p^a_D(t, x, y)\le p^a_{\bH}(t, x, y), \quad a>0,
(t, x, y)\in (0, \infty)\times \bH_b\times \bH_b.
\end{equation}

Let $t_0:=1 \vee b^2$.
 It follows from Theorem \ref{t:main} that we
only need to prove the theorem for $t>3t_0$. Now we suppose
$t>3t_0$. For any $x, y\in D$, we define $x_0$ and $y_0$ as in
\eqref{e:x0y0}.

By the semigroup property and Lemma
\ref{L:rationew},  we have
\begin{eqnarray*}
&&  p^a_D(t,x,y)
 = \int_D \int_D p^a_D(t_0,x,z)p^a_D(t-2t_0,z,w)p^a_D(t_0,w,y)dzdw\\
&&\le  c_{1}\,  (1\wedge \delta_D(x)) (1\wedge
\delta_D(y)) \int_{D \times D }
 h^a_{1/(25C_2)}(t_0, x, z)
  p^a_D
 (t-2t_0,z,w)
 h^a_{1/(25C_2)}(t_0, w, y)  dzdw.
\end{eqnarray*}
It follows from Theorem \ref{t:2sbhalfspace} with $T=1$ and \eqref{e:bdbyhs},
\begin{align*}
&p^a_D(t-2t_0,z,w) \le p^a_{\bH}(t-2t_0,z,w)\\
&\le c_2
 \begin{cases}
 \left(1\wedge
\frac{\delta_{\bH}(z)}{\sqrt{t-2t_0}} \right)\left(1\wedge
\frac{\delta_{\bH}(w)}{\sqrt{t-2t_0}}  \right) \left( (t-2t_0)^{-d/2}
e^{-|z-w|^2/(C_3(t-2t_0))}+\left( \frac{a^\alpha (t-2t_0)}{|z-w|^{d+\alpha}}
\wedge (t-2t_0)^{-d/2} \right)\right),
 \\
 \hskip 4truein ~ \forall (t-2t_0) \in
(0, a^{2\alpha/(\alpha-2)}];
\\
  \left(1\wedge \frac{\delta_{\bH}(z) \wedge
(a^{-1}\delta_{\bH}(z))^{\alpha/2}}{\sqrt{t-2t_0}} \right)\left(1\wedge
\frac{\delta_{\bH}(w)\wedge (a^{-1}\delta_{\bH}(w))^{\alpha/2}}{\sqrt{t-2t_0}}
\right)\left((a^\alpha (t-2t_0))^{-d/\alpha} \wedge \frac{a^\alpha (t-2t_0)
}{|z-w|^{d+\alpha} } \right), \\
 \hskip 4truein ~\forall (t-2t_0) \ge a^{2\alpha/(\alpha-2)} \end{cases}
\end{align*}
where $C_3$ is the constant in Theorem \ref{t:2sbhalfspace} with
$T=1$.
 Put $A=(C_3\vee (25C_2))$ where $C_{2}$ is the constant in
Theorem \ref{t:main} with
$T=t_0$.
  Applying Theorem
\ref{t:2sbhalfspace} with
$T=1$
 again, we get
$$
p^a_D(t-2t_0,z,w)\le c_3 p^a_{\bH}(t-2t_0, A^{-2}z, A^{-2}w)
$$
and so, by Theorem \ref{t:main}
\begin{align*}
 &  p^a_D(t,x,y) \\
 \le & c_{4}\,  (1\wedge \delta_D(x)) (1\wedge
\delta_D(y)) \int_{D \times D }
  h^a_{1/A^4}(t_0, x_0, z)
 p^a_{\bH}(t-2t_0,{A}^{-2}z,{A}^{-2}w)
  h^a_{1/A^4}(t_0, w, y_0) dzdw\nn\\
    \le & c_{5}\,  (1\wedge \delta_D(x)) (1\wedge
\delta_D(y)) \int_{\bH_{-b/2} \times \bH_{-b/2}}
 \left(t_0^{-d/2}e^{-|x_0-z|^2/(A^4t_0)}+\left( \frac{a^{\alpha}t_0}{|x_0-z|^{d+\alpha}}
\wedge  t_0^{-d/2}\right) \right)
 \nn\\
&  \hskip 0.2truein \times p^a_{\bH}(t-2t_0,{A}^{-2}z,{A}^{-2}w)
 \left(t_0^{-d/2}e^{-|w-y_0|^2/(A^4t_0)}+\left( \frac{a^{\alpha}t_0}{|w-y_0|^{d+\alpha}}
\wedge  t_0^{-d/2}\right) \right)  dzdw.\nn
    \end{align*}
Thus, by a change of variable, and
using \eqref{e:bdbyhs} and
Theorem \ref{t:main},
the above is less than or equal to  $(1\wedge \delta_D(x)) (1\wedge
\delta_D(y))$ times
   \begin{align}
    & c_{6}\,   \int_{\bH_{-b/(2A^2)} \times \bH_{-b/(2A^2)}}
 (1\wedge \delta_{\bH_{-b/(2A^2)}}(z))
 (1\wedge \delta_{\bH_{-b/(2A^2)}}(A^{-2}x_0))
 \nn\\
&   \times
\left(t_0^{-d/2}e^{-|x_0-z|^2/(A^4t_0)}+\left( \frac{a^{\alpha}t_0}{|x_0-z|^{d+\alpha}}
\wedge  t_0^{-d/2}\right) \right)
  p^a_{\bH_{-b/(2A^2)}}(t-2t_0,z,w)
  \nn\\
& \times
(1\wedge \delta_{\bH_{-b/(2A^2)}}(A^{-2}y_0))
  (1\wedge \delta_{\bH_{-b/(2A^2)}}(w))
  \left(t_0^{-d/2}e^{-|w-y_0|^2/(A^4t_0)}+\left( \frac{a^{\alpha}t_0}{|w-
  y_0|^{d+\alpha}}
\wedge  t_0^{-d/2}\right) \right)
   dzdw\nn\\
      \le & c_{7}\, \int_{\bH_{-b/(2A^2)} \times \bH_{-b/(2A^2)}}
         p^a_{\bH_{-b/(2A^2)}}(t_0, A^{-2} x_0,     z)
  p^a_{\bH_{-b/(2A^2)}}(t-2t_0,z,w)p^a_{\bH_{-b/(2A^2)}}(t_0, w, A^{-2} y_0)
 dzdw\nn\\
      = &
      c_7\,   p^a_{\bH_{-b/(2A^2)}}(t,A^{-2} x_0, A^{-2} y_0).\nn
 \end{align}
 Now using \eqref{e:scaling} and  Theorem
 \ref{t:2sbhalfspace} with $T=A^{-4}(1\wedge M^{2\alpha/(2-\alpha)})t_0$, we get
\begin{align*}
&
p^a_D(t, x, y)\le
c_8 (1\wedge \delta_D(x))(1\wedge
\delta_D(y))
p^{A^{2(\alpha-2)/\alpha} a}
_{\bH_{-b/2}}(A^{4} t,  x_0,  y_0)\nn\\
 &\le
 c_9
\begin{cases}
(1\wedge \delta_D(x))
\left(\frac{\delta_{\bH_{-b/2}}(x_0)}{\sqrt{t}} \wedge 1 \right)
 (1\wedge \delta_D(y))
\left(\frac{\delta_{\bH_{-b/2}}(y_0)}{\sqrt{t}} \wedge 1 \right) \\
 \quad\times
\left( t^{-d/2} e^{-|x-y|^2/( c_{10}t)} +\left( \frac{a^\alpha
t}{|x-y|^{d+\alpha}} \wedge t^{-d/2} \right)\right) &\hbox{for } t
\in(3t_0,
t_0a^{-2\alpha/(2-\alpha)}],   \\
 (1\wedge \delta_D(x)) \left(\frac{\delta_{\bH_{-b/2}}(x_0)
 \wedge(a^{-1}\delta_{\bH_{-b/2}}(x_0))^{\alpha/2}
}{\sqrt{t}} \wedge 1 \right) (1\wedge \delta_D(y))
 \\ \quad\times
\left(\frac{\delta_{\bH_{-b/2}}(y_0)\wedge(a^{-1}
\delta_{\bH_{-b/2}}(y_0))^{\alpha/2}}{\sqrt{t}} \wedge 1 \right)
\left((a^\alpha t)^{-d/\alpha} \wedge \frac{a^\alpha t
}{|x-y|^{d+\alpha} } \right) &\hbox{for } t>
t_0/a^{2\alpha/(2-\alpha)}
\end{cases}\\
 &\le
 c_{11}
\begin{cases}
(1\wedge \delta_D(x)) \left(\frac{\delta_{{\bH}}(x_0)}{\sqrt{t}}
\wedge 1 \right)
 (1\wedge \delta_D(y))
\left(\frac{\delta_{{\bH}}(y_0)}{\sqrt{t}} \wedge 1 \right)
 \\ \quad\times
 \left( t^{-d/2}
e^{-|x-y|^2/( c_{10}t)}+\left( \frac{a^\alpha t}{|x-y|^{d+\alpha}}
\wedge t^{-d/2} \right)\right) \quad &\hbox{for }   t \in (3t_0,
t_0a^{-2\alpha/(2-\alpha)}];\\
 (1\wedge \delta_D(x)) \left(\frac{\delta_{{\bH}}(x_0)
\wedge(a^{-1}\delta_{\bH}(x_0))^{\alpha/2}}{\sqrt{t}} \wedge 1
\right) (1\wedge \delta_D(y))
 \\ \quad\times
\left(\frac{\delta_{{\bH}}(y_0)\wedge(a^{-1}
\delta_{{\bH}}(y_0))^{\alpha/2}}{\sqrt{t}} \wedge 1 \right)
\left((a^\alpha t)^{-d/\alpha} \wedge \frac{a^\alpha t
}{|x-y|^{d+\alpha} } \right) &\hbox{for } t>
t_0/a^{2\alpha/(2-\alpha)}.
\end{cases}
\end{align*}

In the case when $t> M^{2\alpha/(2-\alpha)}\,
t_0a^{2\alpha/(\alpha-2)}$, since $  M^{2\alpha/(2-\alpha)}\,
t_0a^{2\alpha/(\alpha-2)} \ge t_0$, the desired result follows from
\eqref{e:bdbyhs}, Lemma \ref{l:piw2}, Theorem  \ref{t:2sbhalfspace}
and Remark \ref{R:1.5}(ii).
 In the case when $3t_0<t\le
M^{2\alpha/(2-\alpha)}\, t_0a^{2\alpha/(\alpha-2)}$, the desired
upper bound follows from \eqref{e:bdbyhs},
Theorem \ref{t:2sbhalfspace}, Remark \ref{R:1.5}(ii)
and \cite[Lemma 2.2]{CT} (with $\alpha$ there replaced by 2).

The lower bound can be proved similarly. We omit the details.
\qed

\section{Green function estimates}

In this section, we give the full proof of Theorem
\ref{t:gf-estimates}. Throughout this section, $D$ is a fixed
half-space-like $C^{1,1}$ domain with $C^{1, 1}$ characteristics
$(R_0, \Lambda_0)$ and $\bH_b\subset D\subset \bH$
for some $b>0$ such that the path distance in $D$ is comparable to
the Euclidean distance with characteristic $\lambda_0$.
 We first establish a few lemmas.

Recall that $\phi_a(r)=r\wedge (r/a)^{\alpha/2}$. When $a=1$, we
simply denote $\phi_1$ by $\phi$;
 that is, $\phi(r)= r\wedge r^{\alpha/2}$.

\begin{lemma}\label{l:phi_a}
For every $r\in (0, 1]$ and every open subset $U$ of $\R^d$,
\begin{equation}\label{e:12}
\frac12\left(1\wedge \frac{r^2\phi (\delta_U(x))\phi(\delta_U(y))}{|x-y|^\alpha}\right)
\leq
\left(1\wedge \frac{r \phi (\delta_U(x))}{|x-y|^{\alpha/2}}\right)
\left(1\wedge \frac{r \phi (\delta_U(y))}{|x-y|^{\alpha/2}}\right)
\leq    1\wedge \frac{r^2\phi (\delta_U(x))\phi(\delta_U(y))}{|x-y|^\alpha} .
\end{equation}
\end{lemma}
\pf  The second inequality holds trivially.
Without loss of generality, we assume $\delta_U(x)\leq \delta_U(y)$.
If both $\frac{r \phi (\delta_U(x))}{|x-y|^{\alpha/2}}$
and $\frac{r \phi (\delta_U(y))}{|x-y|^{\alpha/2}}$ are less than 1
or if both are large than one,
$$
\left(1\wedge \frac{r \phi (\delta_U(x))}{|x-y|^{\alpha/2}}\right)
\left(1\wedge \frac{r \phi (\delta_U(y))}{|x-y|^{\alpha/2}}\right)
=  1\wedge \frac{r^2\phi (\delta_U(x))\phi(\delta_U(y))}{|x-y|^\alpha} .
$$
So we only need to consider the case when $\frac{r \phi
(\delta_U(x))}{|x-y|^{\alpha/2}}\leq 1< \frac{r \phi
(\delta_U(y))}{|x-y|^{\alpha/2}}$. Note that
$ \phi (\delta_U (y))\leq \phi (\delta_U(x)+|x-y|)$.
If $\delta_U(x)\geq |x-y|$, then $\phi (\delta_U(y))\leq \phi (2\delta_U(x))
\leq 2 \phi (\delta_U(x))$
and so
$$  1\wedge \frac{r^2\phi (\delta_U(x))\phi(\delta_U(y))}{|x-y|^\alpha}
\leq 1\wedge
2 \left(\frac{r \phi
(\delta_U(x))}{|x-y|^{\alpha/2}} \right)^2
 \leq 2 \left( 1\wedge \frac{r \phi (\delta_U(x))}{|x-y|^{\alpha/2}}\right).
 $$
When $\delta_U(x)<|x-y|$, then $\phi (\delta_U(y))
\leq \phi (2|x-y|)\leq 2  |x-y|^{\alpha/2}$ and so
$$  1\wedge \frac{r^2\phi (\delta_U(x))\phi(\delta_U(y))}{|x-y|^\alpha}
\leq 1\wedge \frac{2 r^2\phi (\delta_U(x)) |x-y|^{\alpha/2} }{|x-y|^\alpha}
\leq 2 \left( 1\wedge \frac{r \phi (\delta_U(x))}{|x-y|^{\alpha/2}}\right)
$$ where the assumption $r \le 1$ is used in the last inequality.
This establishes the first inequality of \eqref{e:12}. \qed

For every open subset $U$ of $\R^d$ and $a >0$, let
\begin{equation}\label{e: q^a_U}
q^a_U(t,x,y):=\left( 1\wedge \frac{\phi_a(\delta_U(x))}{\sqrt{t}}\right)
\left( 1\wedge \frac{\phi_a(\delta_U(y))}{\sqrt{t}}\right)
\left( (a^{\alpha} t)^{-d/\alpha} \wedge \frac{a^{\alpha}t}{|x-y|^{d+\alpha}} \right).
\end{equation}

The following lemma is a direct consequence of (the proof of)
Proposition \ref{C:1.2},
Theorem \ref{t:main_hsl} and Remark \ref{R:1.5}(ii).

\begin{lemma}\label{l:uo1}
For every positive constants $c_1, c_2$, there exists $c_3=c_3(c_1,
c_2)>1$ such that for every $a>0$,
$t \le c_1 a^{-2\alpha/(2-\alpha)}$,
every open subset $U$ of $\R^d$ and $x,y \in U$ with
$|x-y|\ge a^{-\alpha/(2-\alpha)}$,
 \begin{equation}\label{e:1}
 c_3^{-1}\left(1\wedge \frac{\delta_U(x)}{\sqrt{t}}\right)
\left(1\wedge \frac{\delta_U(y)}{\sqrt{t}}\right) h^a_{c_2}(t, x, y)
 \le
 q^a_U (t,x,y) \le
c_3 \left(1\wedge \frac{\delta_U(x)}{\sqrt{t}}\right) \left(1\wedge
\frac{\delta_U(y)}{\sqrt{t}}\right) h^a_{c_2}(t, x, y ).
\end{equation}
Under the assumption of Theorem \ref{t:main}, there is a constant
$c=c(M, R_0, \Lambda_0,
 \lambda_0,  \alpha, b)\geq 1$ such that
$$
c^{-1} q^a_D (t, x, y) \leq p^a_D(t, x, y) \leq c q^a_D(t, x, y)
$$
holds for every $a\in (0, M]$,
$t < \infty$
$x, y \in D$ with $|x-y|\ge a^{-\alpha/(2-\alpha)}$.
\end{lemma}

Observe that
\begin{equation}\label{e:qersa}
\phi_a(\delta_D(\lambda  x))= \big(\lambda \delta_{\lambda^{-1}D}
(x)\big) \wedge \big(\lambda^{\alpha/2}
a^{-\alpha/2}\delta_{\lambda^{-1}D} (x)^{\alpha/2}\big) \qquad
\text{for every }  \lambda>0.
\end{equation}
Let $x_a:=a^{\alpha/(2-\alpha)}x$,  $y_a:=a^{\alpha/(2-\alpha)}y$
and $D_a:= a^{\alpha/(2-\alpha)}D$. By \eqref{e:qersa},
\begin{eqnarray} \label{e:qersa1}
\phi_a(\delta_D(x))=
\phi_a(\delta_D(a^{-\alpha/(2-\alpha)}  x_a))
=a^{-\alpha/(2-\alpha)} \phi(\delta_{D_a}(x_a))
\end{eqnarray}
and so, for every $s>0$,
\begin{eqnarray}\label{e:qersa2}
q^a_D(   a^{-2\alpha/(2-\alpha)}s ,x,y)\,=\,q^a_D(
a^{-2\alpha/(2-\alpha)}s ,a^{-\alpha/(2-\alpha)}
x_a,a^{-\alpha/(2-\alpha)}  y_a)\,=\,a^{\alpha d/(2-\alpha)} q^1_{D_a} (  s ,x_a,y_a).
\end{eqnarray}

We recall that $f^a_D(x,y)$ is defined in \eqref{e:f^a}.

\begin{lemma}\label{l:q^a_U}
For every $d\geq 1$ and $x,y \in D$,
$\int_0^\infty  q^a_D(t,x,y)  dt \asymp f^a_D(x,y)$,
where the implicit constants are independent of $D$.
\end{lemma}

\pf Let $U$ be an arbitrary open subset of $\R^d$. We first consider
the case $a=1$ and prove the lemma for $U$. By a change of variable
  $u= \frac{|x-y|^\alpha}{t}$, we have
\begin{align}
&\int_0^\infty  q^1_U(t,x,y)  dt
\nonumber\\
&=\frac1{|x-y|^{d-\alpha}}  \left(\int_0^1+\int_1^\infty\right)
 \left(u^{(d/\alpha)-2}
\wedge u^{-3} \right)  \left(1\wedge \frac{ {\sqrt u} \phi
(\delta_U(x)) }{ |x-y|^{\alpha/2} }\right) \left(1\wedge \frac{
{\sqrt u} \phi (\delta_U(y))} {|x-y|^{\alpha/2} }\right) du
\nonumber\\
&=:I+II.        \label{e:2}
\end{align}

Note that
\begin{eqnarray}
&&\frac1{2 |x-y|^{d-\alpha}}    \left(1\wedge \frac{ \phi
(\delta_U(x))}{ |x-y|^{\alpha/2} }\right) \left(1\wedge \frac{ \phi
(\delta_U(y))}{ |x-y|^{\alpha/2} }\right)\nonumber\\
&=& \frac1{|x-y|^{d-\alpha}}  \int_{1}^\infty
  u^{-3} \, \left(1\wedge
\frac{   \phi (\delta_U(x))}{ |x-y|^{\alpha/2} }\right)
\left(1\wedge \frac{  \phi (\delta_U(y))}{
|x-y|^{\alpha/2} }\right) du \nonumber\\
&\le& II
\,=\,
 \frac1{|x-y|^{d-\alpha}}  \int_{1}^\infty
  u^{-2} \, \left(u^{-1/2} \wedge
\frac{   \phi (\delta_U(x))}{ |x-y|^{\alpha/2} }\right)
\left(u^{-1/2} \wedge \frac{  \phi (\delta_U(y))}{
|x-y|^{\alpha/2} }\right) du \nonumber\\
&\leq&
 \frac1{|x-y|^{d-\alpha}}
\int_{1}^\infty  u^{-2} \,
  \left(1 \wedge
\frac{   \phi (\delta_U(x))}{ |x-y|^{\alpha/2} }\right) \left(1
\wedge \frac{  \phi (\delta_U(y))}{
|x-y|^{\alpha/2} }\right) du \nonumber\\
&=& \frac1{ |x-y|^{d-\alpha}}    \left(1\wedge \frac{
\phi(\delta_U(x))}{ |x-y|^{\alpha/2} }\right) \left(1\wedge \frac{
\phi (\delta_U(y))}{ |x-y|^{\alpha/2} }\right)  . \label{e:4}
\end{eqnarray}

    \bigskip\noindent
 (i) Assume $d>\alpha$. Observe that
 \begin{eqnarray}
I &\leq & \frac{1} {|x-y|^{d-\alpha}}
 \left( 1\wedge \frac{  \phi (\delta_U(x))}{ |x-y|^{\alpha/2}   }\right)
\left(1\wedge \frac{ \phi (\delta_U(y))}{ |x-y|^{\alpha/2}   }
\right) \int_{0}^1  u^{(d/\alpha)-2} du \nonumber \\
 &\leq &  \frac{\alpha}{d-\alpha} \frac{1} {|x-y|^{d-\alpha}}
 \left(1\wedge \frac{  \phi (\delta_U(x))}{ |x-y|^{\alpha/2}   }\right)
\left( 1\wedge \frac{ \phi (\delta_U(y))}{ |x-y|^{\alpha/2}} \right)
.  \label{e:5}
\end{eqnarray}
So by   \eqref{e:2}--\eqref{e:5},
\begin{equation}\label{e:6}
  \int_0^\infty  q^1_U(t,x,y)  dt\asymp   \frac{1} {|x-y|^{d-\alpha}}
 \left(1\wedge \frac{  \phi (\delta_U(x))}{ |x-y|^{\alpha/2}   }\right)
\left( 1\wedge \frac{ \phi (\delta_U(y))}{ |x-y|^{\alpha/2}} \right)
.
\end{equation}

\bigskip

 For the rest of the proof, we assume without loss of generality that
$\delta_U(x) \le \delta_U(y)$ and define
$$ u_0:=  \frac{
\phi(\delta_U(x)) \phi ( \delta_U(y)) }{ |x-y|^{\alpha} }.
$$

\medskip \noindent
(ii) Now assume $d=\alpha=1$.  We have by Lemma \ref{l:phi_a},
\begin{eqnarray}
I &\asymp&
\int_{0}^{1 } u^{-1} {\bf 1}_{\{u\geq
1/u_0\}} du +  \int_{0}^{1 } u_0 {\bf 1}_{\{u<
1/u_0\}} du \nonumber \\
&=& \log (u_0 \vee 1) + u_0 \left( (1/u_0)\wedge 1 \right)
= \log (u_0 \vee 1) +  (   u_0 \wedge 1  ). \label{e:11}
\end{eqnarray}
Now by Lemma \ref{l:phi_a}, \eqref{e:2}-\eqref{e:4} and \eqref{e:11}, we
have
   \begin{align*}
        \int_0^\infty  q^1_U(t,x,y)  dt
        \asymp \log(u_0\vee 1) + 1\wedge u_0
        \asymp \log(1 + u_0).
    \end{align*}

\medskip\noindent
(iii) Lastly we consider the case $d=1<\alpha<2$.
 By Lemma \ref{l:phi_a},
 \begin{eqnarray*}
I &\asymp&
\frac{1}{|x-y|^{1-\alpha}} \left(
\int_{0}^{1 } u^{(1/\alpha)-2}  {\bf
1}_{\{u\geq 1/u_0\}} du +  \int_{0}^{1 } u_0
u^{(1/\alpha)-1} {\bf 1}_{\{u<
1/u_0\}} du \right)
\\
 &=&  \frac{1}{|x-y|^{1-\alpha}}\left( \frac{\alpha}{\alpha-1}
 \left( (u_0\vee 1)^{1-(1/\alpha)}-1\right) + \alpha u_0
 (u_0\vee 1)^{-1/\alpha}
 \right).
\end{eqnarray*}
 Hence by  \eqref{e:2}-\eqref{e:4}, Lemma \ref{l:phi_a}
    and the last display we have
    \begin{align*}
    &\int_0^\infty  q^1_U(t,x,y)  dt \\
  &\asymp \frac{1}{|x-y|^{1-\alpha}}(1\wedge u_0) + \frac{1}{|x-y|^{1-\alpha}}
  \left( \left( (u_0\vee 1)^{1 - (1/\alpha)} - 1 \right) + u_0 (u_0\vee 1)^{-1/\alpha}  \right) \\
  &\asymp \frac{1}{|x-y|^{1-\alpha}}\left( u_0\wedge u_0^{1-(1/\alpha)} \right)
 = \frac{
\phi(\delta_U(x)) \phi ( \delta_U(y)) }{ |x-y|  } \wedge \left(
\phi(\delta_U(x))  \phi (\delta_U(y))\right)^{(\alpha-1)/\alpha}.
    \end{align*}
Thus we have proved the lemma for any open set $U$ and $a=1$. For
general $a>0$,  we have by \eqref{e:qersa1} and \eqref{e:qersa2}
that
       \begin{align*}
&   \int_0^\infty
     q^a_D(  t ,x,y) dt
     \,=\,
       a^{-2\alpha/(2-\alpha)}\int_0^\infty
     q^a_D(   a^{-2\alpha/(2-\alpha)}s ,x,y) ds
  \,=\, a^{\alpha(d-2)/(2-\alpha)}  \int_0^\infty   q^1_{D_a} (  s ,x_a,y_a) ds\nn\\
  \asymp&  a^{\alpha(d-2)/(2-\alpha)}  \begin{cases}
 \frac{1} {|x_a-y_a|^{d-\alpha}}
 \left(1\wedge \frac{  \phi (\delta_{D_a}(x_a))}{ |x_a-y_a|^{\alpha/2}   }\right)
\left( 1\wedge \frac{ \phi (\delta_{D_a}(y_a))}{
|x_a-y_a|^{\alpha/2}}
\right)  &\hbox{when } d>\alpha, \smallskip  \\
\log \left( 1+ \frac{ \phi(\delta_{D_a}(x_a)) \phi (
\delta_{D_a}(y_a)) }{ |x_a-y_a|^{\alpha} }\right)
  &\hbox{when } d=1=\alpha,  \smallskip \\
\frac{ \phi(\delta_{D_a}(x_a)) \phi ( \delta_{D_a}(y_a)) }{
|x_a-y_a|  } \wedge \left( \phi(\delta_{D_a}(x_a))  \phi (
\delta_{D_a}(y_a))\right)^{(\alpha-1)/\alpha}
 &\hbox{when } d=1<\alpha .
\end{cases}
\\
  =&  a^{\alpha(d-2)/(2-\alpha)}  \begin{cases}
\frac{a^{-(d-\alpha)\alpha/(2-\alpha)}  } {|x-y|^{d-\alpha}}
\left(1\wedge \frac{
a^{\alpha/(2-\alpha)}\phi_a(\delta_D(x))}{a^{\alpha^2/2(2-\alpha)}
|x-y|^{\alpha/2}   }\right) \left( 1\wedge \frac{
a^{\alpha/(2-\alpha)}\phi_a(\delta_D(y))}{
a^{\alpha^2/2(2-\alpha)}|x-y|^{\alpha/2}}
\right)  &\hbox{when } d>\alpha, \smallskip  \\
\log \left( 1+ \frac{
a^{2}\phi_a(\delta_D(x)) \phi_a ( \delta_{D}(y)) }{a |x-y| }\right)
  &\hbox{when } d=1=\alpha,  \smallskip \\
\frac{ a^{2\alpha/(2-\alpha)}\phi_a(\delta_D(x)) \phi_a (
\delta_{D}(y)) }{a^{\alpha/(2-\alpha)} |x-y|  } \wedge \left(
a^{2\alpha/(2-\alpha)}\phi_a(\delta_D(x))  \phi_a (
\delta_D(y))\right)^{(\alpha-1)/\alpha}
 &\hbox{when } d=1<\alpha
\end{cases}
\\
=&f^a_D(x,y).
\end{align*}
\qed

\begin{lemma}\label{l:I_est}
For every $c>0$,  when $|x-y|  \le  a^{-\alpha/(2-\alpha)}$,
\begin{eqnarray*}
&&\int_0^{a^{-2\alpha/(2-\alpha)}}
 \left(1\wedge \frac{\delta_D(x)}{\sqrt{t}}\right)
\left(1\wedge \frac{\delta_D(y)}{\sqrt{t}}\right) \left[t^{-d/2}
 e^{-c\frac{|x-y|^2}{t}}+\left(\frac{ a^\alpha t}{|x-y|^{d +
\alpha}}\wedge
t^{-d/2}\right)\right]\, dt\\
&\asymp&
\begin{cases} |x-y|^{2-d} \left(1\wedge \frac{
\delta_D(x) \delta_D(y)}{ |x-y|^2 }\right) &\hbox{when } d \ge 3 , \smallskip  \\
\log (1+ \frac{a^{2\alpha/(\alpha-2)}\wedge (\delta_D(x)
\delta_D(y))}{ |x-y|^2 } ) &\hbox{when }
d=2, \smallskip \\
a^{\alpha/(\alpha-2)}\wedge\left( \delta_D(x) \delta_D
(y)\right)^{1/2} \, \wedge \, \frac{ \delta_D(x) \delta_D (y)}{
|x-y| } &\hbox{when } d=1,
\end{cases}
\end{eqnarray*}
where the implicit constant depend only on $c, \alpha$ and $d$.
\end{lemma}

\pf We first consider the case $a=1$ and assume $U$
is an arbitrary open set and $x,y \in U$ with $|x-y| \le 1$.
Using the change of variables $u=\frac{|x-y|^2}{t}$, we have
\begin{align*}
&\int_0^1 \left(1\wedge \frac{\delta_U(x)}{\sqrt{t}}\right)
\left(1\wedge \frac{\delta_U(y)}{\sqrt{t}}\right) \left[t^{-d/2}
e^{-c_1\frac{|x-y|^2}{t}}+\left(\frac{t}{|x-y|^{d+\alpha}}\wedge
t^{-d/2}\right)\right]\, dt\\
=& |x-y|^{2-d} \, \left(\int_{|x-y|^2}^{2}+\int_2^\infty\right) \left(1\wedge
\frac{\sqrt{u}\delta_U(x)}{|x-y|}\right) \left(1\wedge
\frac{\sqrt{u}\delta_U(y)}{|x-y|}\right) \left[ u^{d/2} e^{-c_1u}+
\left(\frac{|x-y|^{2-\alpha}}{u}\wedge u^{d/2}\right) \right]
\frac{du}{u^2}\\
=:&I_1+I_2\, .
\end{align*}
Note that since $|x-y|^{2-\alpha}\le 1$,
 for $u\geq 2$,
$\frac{|x-y|^{2-\alpha}}{u}\wedge
u^{d/2}=\frac{|x-y|^{2-\alpha}}{u}$.
Thus
for any $d \ge 1$,
\begin{eqnarray*}
I_2&=& |x-y|^{2-d}\int_2^{\infty} \left(u^{-1/2}\wedge
\frac{\delta_U(x)}{|x-y|}\right) \left(u^{-1/2}\wedge
\frac{\delta_U(y)}{|x-y|}\right) \left[u^{d/2} e^{-c_1u}+
\frac{|x-y|^{2-\alpha}}{u}\right]\, \frac{du}{u}\\
&\le & |x-y|^{2-d}\left(1\wedge \frac{\delta_U(x)}{|x-y|} \right)
\left(1\wedge \frac{\delta_U(y)}{|x-y|}\right)
\int_2^{\infty}\left(u^{d/2-1} e^{-c_1u}+u^{-2}\right)\, du\\
&\le & c_2 |x-y|^{2-d}\left(1\wedge \frac{\delta_U(x)}{|x-y|}
\right) \left(1\wedge \frac{\delta_U(y)}{|x-y|}\right) \,
\end{eqnarray*}
and
\begin{eqnarray*}
I_2&\ge& |x-y|^{2-d}\int_2^{\infty} \left(1\wedge
\frac{\delta_U(x)}{|x-y|}\right) \left(1\wedge \frac{\delta_U(y)}
{|x-y|}\right)\left[u^{d/2} e^{-c_1u}+\frac{|x-y|^{2-\alpha}}
{u}\right]\, \frac{du}{u^2}\\
&\ge &|x-y|^{2-d}\left(1\wedge \frac{\delta_U(x)}{|x-y|}\right)
\left(1\wedge \frac{\delta_U(y)}{|x-y|}\right) \int_2^{\infty}
u^{d/2-2} e^{-c_1u}\, du\\
&\ge  &c_3 |x-y|^{2-d}\left(1\wedge \frac{\delta_U(x)}{|x-y|}
\right) \left(1\wedge \frac{\delta_U(y)}{|x-y|}\right).
\end{eqnarray*}

One the other hand,  since $|x-y|^{2-\alpha}\le 1$, if $u\le 2$,
then
$$
u^{-2} \left[ u^{d/2} e^{-c_1u}+
\left(\frac{|x-y|^{2-\alpha}}{u}\wedge u^{d/2}\right) \right] \asymp
u^{d/2-2}.
$$
Using this and the fact that for every $r\in (0, 2]$,
\begin{equation}\label{e:ch5}
\left(1\wedge \frac{r \delta_U(x)} {|x-y|}\right)\, \left( 1\wedge
\frac{r \delta_U(y)} {|x-y|} \right) \, \le \, 1\wedge
\frac{r^2\delta_U(x)\delta_U(y)} {|x-y|^2} \,\le\,
4 \left(
1\wedge \frac{r \delta_U(x)} {|x-y|} \right)\, \left( 1\wedge
\frac{r \delta_U(y)} {|x-y|} \right),
\end{equation}
we have
$$
I_1 \asymp |x-y|^{2-d} \, \int_{|x-y|^2}^{2} \left(1\wedge
\frac{u\delta_U(x)\delta_U(y)}{|x-y|^2}\right)  u^{d/2-2}\, {du}.
$$

 Let $ u_0:=  \frac{
\delta_U(x) \delta_U(y)}{ |x-y|^2 }.$

\noindent
 (i) When $d \ge 3$,  it is easy to see that
$
I_1\le   |x-y|^{2-d} \left(1\wedge u_0\right) .
$

\noindent
 (ii) Assume $d =2$.
We deal with three cases separately.

(a) $u_0\le 1$: In this case, since $|x-y|\leq 1$, we have
$\delta_U(x)\delta_U(y)\leq 1$ and
$ I_1 \asymp  \int_{|x-y|^2}^{2}u_0du\asymp u_0\asymp \ln(1+u_0)$.

(b) $u_0> 1$ and $|x-y|^2\le 1/u_0$:  In this case we have
$\delta_U(x)\delta_U(y)\leq 1$ and
\begin{align*}
I_1 \asymp&
\int^{u^{-1}_0}_{|x-y|^2}u_0du+\int^2_{u^{-1}_0}u^{-1}du
=u_0(u^{-1}_0-|x-y|^2)+\ln 2 +\ln u_0\\
=&(1-u_0|x-y|^2)+\ln 2 +\ln u_0\asymp \ln(1+u_0).
\end{align*}

(c) $u_0> 1$ and $|x-y|^2> 1/u_0$:  In this case we have
$\delta_U(x)\delta_U(y)\geq 1$ and
\begin{eqnarray*}
I_1 \asymp\int^2_{|x-y|^2}u^{-1}du=\ln 2+ \ln |x-y|^{-2}
\asymp\ln (1+ |x-y|^{-2})=\ln \left(
1+\frac{1\wedge(\delta_U(x)\delta_U(y))}{|x-y|^2} \right).
\end{eqnarray*}

 \noindent
(iii) Now we consider the case $d=1$.
We again deal with three cases separately.

(a) $u_0\le 1$. In this case we have
$$
I_1 \asymp  |x-y|\int_{|x-y|^2}^{2}u_0u^{-1/2}du\asymp
|x-y|u_0(\sqrt 2-|x-y|)\asymp |x-y|u_0.
$$

(b) $u_0> 1$ and $|x-y|^2\le 1/u_0$. In this case we have
\begin{align*}
I_1 \asymp& |x-y|
\int^{u^{-1}_0}_{|x-y|^2}u_0u^{-1/2}du + |x-y|\int^2_{u^{-1}_0}u^{-3/2}du\\
\asymp&u_0|x-y|(u^{-1/2}_0-|x-y|)+ |x-y|(u_0^{1/2}-2^{-1/2})\asymp|x-y|u_0^{1/2}.
\end{align*}

(c) $u_0> 1$ and $|x-y|^2> 1/u_0$. In this case we have
\begin{align*}
I_1 \asymp|x-y|\int^2_{|x-y|^2}u^{-3/2}du \asymp
|x-y|(|x-y|^{-1}-2^{-1/2})\asymp 1-2^{-1/2}|x-y|\asymp 1.
\end{align*}
So we have
\begin{equation}\label{e:dssw1}
I_1+I_2 \asymp
\begin{cases} |x-y|^{2-d} \left(1\wedge \frac{
\delta_U(x) \delta_U(y)}{ |x-y|^2 }\right)        &\hbox{when } d \ge 3 , \smallskip  \\
\log (1+ \frac{1\wedge( \delta_U(x) \delta_U(y))}{ |x-y|^2 } )
&\hbox{when }
d=2,\smallskip \\
1\wedge\left( \delta_U(x) \delta_U (y)\right)^{1/2} \, \wedge \,
\frac{ \delta_U(x) \delta_U (y)}{ |x-y| } &\hbox{when } d=1.
\end{cases}
\end{equation}

Thus we have proved the lemma for any  open set $U$ and $a=1$. For
general $a>0$, we have by \eqref{e:qersa1}, \eqref{e:qersa2} and
\eqref{e:dssw1},
\begin{align}
&   \int_0^{a^{-2\alpha/(2-\alpha)}} \left(1\wedge
\frac{\delta_D(x)}{\sqrt{t}}\right) \left(1\wedge
\frac{\delta_D(y)}{\sqrt{t}}\right) \left[t^{-d/2}
e^{-c_1\frac{|x-y|^2}{t}}+\left(\frac{ a^\alpha
t}{|x-y|^{d+\alpha}}\wedge
t^{-d/2}\right)\right]\, dt    \nn\\
=&    a^{-2\alpha/(2-\alpha)}\int_0^{1} \left(1\wedge
\frac{\delta_D(x)}{a^{-\alpha/(2-\alpha)}\sqrt{s}}\right)
\left(1\wedge \frac{\delta_D(y)}{a^{-\alpha/(2-\alpha)}\sqrt{s}}\right) \nn\\
&\times\left[(a^{-2\alpha/(2-\alpha)}s)^{-d/2}
e^{-c_1\frac{|x-y|^2}{a^{-2\alpha/(2-\alpha)}s}}+ \left(\frac{
a^\alpha a^{-2\alpha/(2-\alpha)}s}{|x-y|^{d+\alpha}}\wedge
(a^{-2\alpha/(2-\alpha)}s)^{-d/2}\right)\right]\, ds
\nn\\
=&    a^{\alpha(d-2)/(2-\alpha)}\int_0^{1} \left(1\wedge
\frac{\delta_{D_a}(x_a)}{\sqrt{s}}\right) \left(1\wedge
\frac{\delta_{D_a}(y_a)}{\sqrt{s}}\right) \left[s^{-d/2}
e^{-c_1\frac{|x_a-y_a|^2}{s}}+ \left(\frac{s}{|x_a-y_a|^{d+\alpha}}
\wedge s^{-d/2}\right)\right]\, ds\nn\\
\asymp&  a^{\alpha(d-2)/(2-\alpha)}  \begin{cases} |x_a-y_a|^{2-d}
\left(1\wedge \frac{ \delta_{D_a}(x_a) \delta_{D_a}(y_a)}{
|x_a-y_a|^2 }\right)
&\hbox{when } d \ge 3 , \smallskip  \\
\log (1+ \frac{1\wedge( \delta_{D_a}(x_a)
\delta_{D_a}(y_a))}{|x_a-y_a|^2 } )
&\hbox{when } d=2,\smallskip \\
1\wedge \left( \delta_{D_a}(x_a) \delta_{D_a}(y_a)\right)^{1/2} \,
\wedge \, \frac{  \delta_{D_a}(x_a) \delta_{D_a}(y_a)}{ |x_a-y_a| }
&\hbox{when } d=1
\end{cases}\nn\\
=&   \begin{cases} |x-y|^{2-d} \left(1\wedge \frac{
\delta_D(x) \delta_D(y)}{ |x-y|^2 }\right) &\hbox{when } d \ge 3 , \smallskip  \\
\log (1+ \frac{a^{2\alpha/(\alpha-2)}\wedge (\delta_D(x)
\delta_D(y))}{ |x-y|^2 } ) &\hbox{when }
d=2, \smallskip \\
a^{\alpha/(\alpha-2)}\wedge\left( \delta_D(x) \delta_D
(y)\right)^{1/2} \, \wedge \, \frac{ \delta_D(x) \delta_D (y)}{
|x-y| } &\hbox{when } d=1.
\end{cases} \nn
\end{align}
\qed

\begin{lemma}\label{l:J_est_d2}
For every $d \ge 2$, there exists $c=c(\alpha, d)>1$ such that, for every $a>0$,
when $|x-y|  \le  a^{-\alpha/(2-\alpha)}$,
\begin{eqnarray*}
\int^\infty_{a^{-2\alpha/(2-\alpha)}} q_D^a(t,x,y)\, dt \le
 c\left(1\wedge \frac{\delta_D(x)\delta_D(y)}{|x-y|^2}
\right).
 \end{eqnarray*}
\end{lemma}

\pf We first consider the case $a=1$ and assume $U$ is an arbitrary
open set and $x,y \in U$ with $|x-y| \le 1$. Let $J:=\int^\infty_{1}
q_U^1(t,x,y)\, dt$.
By a
change of variables $u=\frac{|x-y|^{\alpha}}{t}$,
\begin{align}
J =
|x-y|^{\alpha-d} \int_0^{|x-y|^{\alpha}} \left(1\wedge
\frac{\sqrt{u}(\delta_U(x)\wedge
\delta_U(x)^{\alpha/2})}{|x-y|^{\alpha/2}}\right) \left(1\wedge
\frac{\sqrt{u}(\delta_U(y)\wedge
\delta_U(y)^{\alpha/2})}{|x-y|^{\alpha/2}}\right)
\left(u^{d/\alpha}\wedge u^{-1}\right)\, \frac{du}{u^2}\, . \label{e:J}
\end{align}
Since $|x-y|\le 1$, for $u\in [0, \, |x-y|^\alpha]$,
$u^{d/\alpha}\wedge u^{-1}=u^{d/\alpha}$. Hence
\begin{eqnarray*}
J&\le & |x-y|^{\alpha-d} \left(1\wedge
\frac{\delta_U(x)\wedge \delta_U(x)^{\alpha/2}}{|x-y|^{\alpha/2}}
\right) \left(1\wedge \frac{\delta_U(y)\wedge
\delta_U(y)^{\alpha/2}}
{|x-y|^{\alpha/2}}\right)\int_0^{|x-y|^{\alpha}} u^{d/\alpha-2}\, du\\
&=&
c_1 \left(1\wedge \frac{\delta_U(x)\wedge \delta_U(x)^{\alpha/2}}
{|x-y|^{\alpha/2}}\right) \left(1\wedge \frac{\delta_U(y)\wedge
\delta_U(y)^{\alpha/2}}{|x-y|^{\alpha/2}}\right)\, .
\end{eqnarray*}
Since $|x-y|\le |x-y|^{\alpha/2}\le 1$, we have that
$\frac{1}{|x-y|^{\alpha/2}}\le \frac{1}{|x-y|}$ and so
$
1\wedge \frac{\delta_U(x)\wedge \delta_U(x)^{\alpha/2}}
{|x-y|^{\alpha/2}}\le 1\wedge \frac{\delta_U(x)}{|x-y|}\,$.
Consequently,
  \begin{eqnarray}
J\le
c_1 \left(1\wedge \frac{\delta_U(x)}{|x-y|}\right)
\left(1\wedge \frac{\delta_U(y)}{|x-y|}\right) \le
 2
c_1
\left(1\wedge \frac{\delta_U(x)\delta_U(y)}{|x-y|^2}\right)   \, . \label{e:J11}
  \end{eqnarray}

Thus we have proved the lemma for any  open set $U$ and $a=1$. For
general $a>0$, by \eqref{e:qersa1}, \eqref{e:qersa2} and
\eqref{e:J11}, we have
\begin{align*}
 & \int_{a^{-2\alpha/(2-\alpha)}}^\infty
     q^a_D(  t ,x,y) dt
    = a^{\alpha(d-2)/(2-\alpha)}  \int_1^\infty   q^1_{D_a} (  s ,x_a,y_a) ds
 \\& \le   2
c_1 a^{\alpha(d-2)/(2-\alpha)}
\left(1\wedge \frac{\delta_{D_a}(x_a)\delta_{D_a}(y_a)}{|x_a-y_a|^2}\right)=     2
c_1 \left(1\wedge
\frac{\delta_D(x)\delta_D(y)}{|x-y|^2} \right).
\end{align*}
\qed

\begin{lemma}\label{l:IJ_estd1}
 For every $c>0$, when $d=1$ and $|x-y|  \le  a^{-\alpha/(2-\alpha)}$,
\begin{eqnarray*}
&&\int_0^{a^{-2\alpha/(2-\alpha)}}
 \left(1\wedge \frac{\delta_D(x)}{\sqrt{t}}\right)
\left(1\wedge \frac{\delta_D(y)}{\sqrt{t}}\right) \left( t^{-d/2}
 e^{-c\frac{|x-y|^2}{t}}+\left(\frac{ a^\alpha t}{|x-y|^{d+\alpha}}\wedge
t^{-d/2}\right)\right) \, dt\\
&&+\int^\infty_{a^{-2\alpha/(2-\alpha)}}
q_D^a(t,x,y)\, dt\,\asymp\,
g_D^a(x, y)
\end{eqnarray*}
where the implicit constant depend only on
$c$ and $\alpha$.
\end{lemma}

\pf We first consider the case $a=1$ and assume $U$ is an arbitrary
open set and $x,y \in U$ with $|x-y| \le 1$. Let $J:=\int^\infty_{1}
q_U^1(t,x,y)\, dt$ and
$$I:=\int_0^{1}
 \left(1\wedge \frac{\delta_U(x)}{\sqrt{t}}\right)
\left(1\wedge \frac{\delta_U(y)}{\sqrt{t}}\right)
\left( t^{-1/2}
e^{-c_1\frac{|x-y|^2}{t}}+\left(\frac{  t}{|x-y|^{1+\alpha}}\wedge
t^{-1/2}\right)\right) \, dt.
$$
By Lemma \ref{l:I_est}, $I \asymp  1\wedge \left( \delta_D(x)
\delta_D (y)\right)^{1/2} \, \wedge \, \frac{  \delta_D(x) \delta_D
(y)}{ |x-y| }.$ Using Lemma \ref{l:phi_a} and \eqref{e:J},  we get
that
$$
\int_1^\infty
q_U^1(t,x,y)\, dt \asymp |x-y|^{\alpha-1}
\int_0^{|x-y|^{\alpha}} \left(1\wedge\frac{u\phi(\delta_U(x))
\phi(\delta_U(y))}{|x-y|^{\alpha}}\right) u^{1/\alpha-2}\, du.
$$
Put
$u_0:=\frac{\phi(\delta_U(x)) \phi(\delta_U(y))}{|x-y|^{\alpha}}$.
Then we have
\begin{eqnarray*}
J&\asymp&|x-y|^{\alpha-1}\left(
u_0
\int_0^{|x-y|^\alpha\wedge
u_0^{-1}}u^{1/\alpha-1}\, du +\int_{|x-y|^\alpha\wedge
u_0^{-1}}^{|x-y|^\alpha}u^{1/\alpha-2}\, du\right).
\end{eqnarray*}

Without loss of generality, we assume $\delta_U(x) \le \delta_U(y)$.
Note that, since $|x-y| \le 1$, if $\delta_U(x) \le 1$ then
$\delta_U(y) \le 2$, and if $\delta_U(x) > 1$ then $1 <\delta_U(x)
\le \delta_U(y) \le 2 \delta_U(x)$ and $\delta_U(x)\delta_U(y) \ge
|x-y|^2$.

Now we look at three separate cases.

\noindent
(i) $\alpha\in (1, 2)$: In this case we have
\begin{eqnarray*}
J&\asymp&|x-y|^{\alpha-1}\left(\alpha
u_0 \left(|x-y|\wedge u^{-1/\alpha}_0 \right)
 + \frac{\alpha}{\alpha-1}
\left(|x-y|^\alpha\wedge u_0^{-1}\right)^{(1-\alpha)/\alpha}
-\frac{\alpha}{\alpha-1}|x-y|^{1-\alpha}
\right)\\
&\asymp&\phi(\delta_U(x))\phi(\delta_U(y))  \wedge
(\phi(\delta_U(x))\phi(\delta_U(y)))^{(\alpha-1)/\alpha}.
\end{eqnarray*}
Thus
\begin{eqnarray*}
I+J&\asymp&
\begin{cases}
\left( \delta_U(x) \delta_U (y)\right)^{1/2}
&\hbox{when } \delta_U(x) \le 1,   \delta_U(x)\delta_U(y) \ge |x-y|^2,
\smallskip  \\
\frac{  \delta_U(x) \delta_U (y)}{ |x-y| }
&\hbox{when } \delta_U(x) \le 1,  \delta_U(x)\delta_U(y) \le |x-y|^2, \smallskip \\
\left( \delta_U(x) \delta_U (y)\right)^{(\alpha-1)/2}
 &\hbox{when } \delta_U(x) > 1 \end{cases}\\
&=& \left( \delta_U(x) \delta_U (y)\right)^{1/2}
   \wedge  \left( \delta_U(x) \delta_U (y)\right)^{(\alpha-1)/2}
   \wedge \frac{  \delta_U(x) \delta_U (y)}{ |x-y| } .
\end{eqnarray*}

\noindent
(ii) $\alpha=1$: In this case we have
\begin{eqnarray*}
J&\asymp&
\left( u_0 (|x-y|\wedge u^{-1}_0) +\log
\frac{|x-y|^{\alpha}}{|x-y|^\alpha\wedge u_0^{-1}}
\right)\\
&\asymp&\phi(\delta_U(x))\phi(\delta_U(y)) \wedge 1+
\log\left(1 \vee \phi(\delta_U(x))\phi(\delta_U(y))\right)
\asymp\log\left(1+\phi(\delta_U(x))\phi(\delta_U(y))\right) .
\end{eqnarray*}
Thus
\begin{eqnarray*}
I+J&\asymp&
\begin{cases}
\left( \delta_U(x) \delta_U (y)\right)^{1/2}
&\hbox{when } \delta_U(x) \le 1,     \delta_U(x)\delta_U(y) \ge |x-y|^2,
\smallskip  \\
\frac{  \delta_U(x) \delta_U (y)}{ |x-y| }
&\hbox{when } \delta_U(x) \le 1,  \delta_U(x)\delta_U(y) \le |x-y|^2,
\smallskip \\
\log\left(1+\delta_U(x)\delta_U(y)\right)
&\hbox{when } \delta_U(x) > 1 \end{cases}\\
  &\asymp&  \frac{
\delta_U(x) \delta_U (y)}{ |x-y| } \wedge
  \log\left(1+\left( \delta_U(x)\delta_U(y)\right)^{1/2} \right) .
\end{eqnarray*}

\noindent (iii) $\alpha\in (0, 1)$: In this case (note that
$1-1/\alpha$ is negative) we have
\begin{eqnarray*}
J&\asymp&|x-y|^{\alpha-1}\left(\alpha
u_0(|x-y|\wedge u^{-1/\alpha}_0)
+\frac{\alpha}{1-\alpha}|x-y|^{1-\alpha}
-\frac{\alpha}{1-\alpha}(|x-y|^\alpha \wedge
u_0^{-1})^{(1-\alpha)/\alpha} \right)\\
&\asymp&\phi(\delta_U(x))\phi(\delta_U(y)) \wedge 1.
\end{eqnarray*}
Thus
\begin{eqnarray*}
I+J&\asymp&
\begin{cases}
\left( \delta_U(x) \delta_U (y)\right)^{1/2}
\quad &\hbox{when } \delta_U(x) \le 1, \delta_U(x)\delta_U(y) \ge |x-y|^2,
\smallskip  \\
\frac{  \delta_U(x) \delta_U (y)}{|x-y| }
&\hbox{when } \delta_U(x) \le 1,  \delta_U(x)\delta_U(y) \le |x-y|^2,
\smallskip \\
 1
&\hbox{when } \delta_U(x) > 1 \end{cases}\\
&=&  \left( \delta_U(x) \delta_U (y)\right)^{1/2} \wedge \frac{
\delta_U(x) \delta_U (y)}{ |x-y| } \wedge 1.
\end{eqnarray*}

Therefore we have proved the lemma for any arbitrary open set $U$
and $a=1$. The general case $a>0$ now follows from the same scaling
arguments as in the proofs for
 Lemmas \ref{l:q^a_U} and \ref{l:I_est}.
 \qed

\noindent {\bf Proof of  Theorem \ref{t:gf-estimates}}.
Without loss of generality, we assume $M=b=1$. Estimates
\eqref{e:gf-estimates1}
follow from  Theorem \ref{t:main_hsl}, Remark \ref{R:1.5}(ii)
 and Lemmas \ref{l:I_est}--\ref{l:IJ_estd1}.
 Estimates
\eqref{e:gf-estimates2}
follow from Theorem \ref{t:main_hsl} and Lemmas \ref{l:uo1} and \ref{l:q^a_U}.
\qed

\medskip
\noindent {\bf Acknowledgment:} While working on the paper
\cite{KSV2}, Z. Vondra\v{c}ek
 obtained
 the Green function
estimates of $p^1_{\bH}$ in the case $d\ge 3$ using Theorem
\ref{t:main_hsl} above. Some of his calculations are incorporated in
the proofs of Lemmas \ref{l:I_est}--\ref{l:J_est_d2}.

\vskip 0.3truein

{\bf Zhen-Qing Chen}

Department of Mathematics, University of Washington, Seattle,
WA 98195, USA

E-mail: \texttt{zchen@math.washington.edu}

\bigskip

{\bf Panki Kim}

Department of Mathematical Sciences and Research Institute of Mathematics,

Seoul National University,
San56-1 Shinrim-dong Kwanak-gu,
Seoul 151-747, Republic of Korea

E-mail: \texttt{pkim@snu.ac.kr}

\bigskip

{\bf Renming Song}

Department of Mathematics, University of Illinois, Urbana, IL 61801, USA

E-mail: \texttt{rsong@math.uiuc.edu}

\small

\begin{thebibliography}{99}

\bibitem{Ber} J. Bertoin: {\em L\'evy Processes}. Cambridge University
Press, Cambridge (1996).

\bibitem{BBC} K. Bogdan, K. Burdzy and Z.-Q. Chen:
Censored stable processes. {\it Probab. Theory Relat. Fields} {\bf
127} (2003), 89--152.

\bibitem{BGR} K. Bogdan, T. Grzywny and M. Ryznar:
Heat kernel estimates for the fractional Laplacian with Dirichlet
Conditions. {\it  Ann. Probab.}
{\bf 38(5)} (2010), 1901--1923.

\bibitem{CKS}  Z.-Q. Chen, P. Kim, and R. Song:
Heat kernel estimates for Dirichlet fractional Laplacian.
{\em J. European Math. Soc.} {\bf 12} (2010), 1307--1329.


\bibitem{CKS5} Z.-Q. Chen, P. Kim and R. Song:
Heat kernel estimate for $\Delta+\Delta^{\alpha/2}$ in $C^{1,1}$
open sets. To appear in {\it J. London Math. Soc.}

\bibitem{CKS4}  Z.-Q. Chen, P. Kim, and R. Song:
Global heat kernel estimates for relativistic stable processes in
half-space-like open sets. To apeear in {\it Potential Anal.}

\bibitem{CKSV} Z.-Q. Chen, P. Kim, R. Song and Z. Vondra{\v{c}}ek:
Boundary Harnack principle for $\Delta + \Delta^{\alpha/2}$. To
appear in {\it Trans. Amer. Math. Soc.}

\bibitem{CK} Z.-Q. Chen and T. Kumagai:
Heat kernel estimates for stable-like processes on $d$-sets. {\em
Stoch. Proc. Appl.}   {\bf 108} (2003), 27--62.

\bibitem{CK2}  Z.-Q. Chen and  T. Kumagai:
Heat kernel estimates for jump processes of mixed types on metric
measure spaces. {\it Probab. Theory Relat. Fields}  {\bf 140}
(2008), 277--317.

\bibitem{CK08}  Z.-Q. Chen and  T. Kumagai:
A priori H\"older estimate, parabolic Harnack principle and heat
kernel estimates for diffusions with jumps.   {\it Rev. Mat.
Iberoam. \bf 26} (2010), 551-589.

\bibitem{CT}
Z.-Q.~Chen and J. Tokle:  Global heat kernel estimates for
fractional Laplacians in unbounded  open sets. {\em Probab. Theory
Relat. Fields}, DOI 10.1007/s00440-009-0256-0 (online first).

\bibitem{G}  Q.-Y. Guan: Boundary Harnack inequality for regional
fractional Laplacian. \ arXiv:0705.1614v3 [math.PR]

\bibitem{KS1} P. Kim and   R. Song: Potential theory of truncated stable processes.
{\it Math. Z.} {\bf 256} (2007), 139--173.

\bibitem{KSV} P. Kim, R. Song and Z. Vondra\v{c}ek: Boundary Harnack
principle for subordinate Brownian motion. \ {\em Stoch. Proc. Appl.}
{\bf 119} (2009), 1601--1631.


\bibitem{KSV2} P. Kim, R. Song and Z. Vondra\v{c}ek:
Minimal thinness for subordinate Brownian motion in half
space. Preprint.


\bibitem{KSV09} P. Kim, R. Song and Z. Vondra\v{c}ek:
On the potential theory of one-dimensional subordinate Brownian
motions with continuous components.   \ {\it Potential Anal.
\bf 33} (2010), 153-173 2009.


\bibitem{SSV} R. L. Schilling, R. Song and Z. Vondra{\v{c}}ek:
{\em Bernstein Functions: Theory and Applications}. de Gruyter
Studies in Mathematics 37. Berlin: Walter de Gruyter, 2010.

\bibitem{Sil} M. L. Silverstein: Classification of coharmonic and coinvariant functions for a
L\'evy process. {\it Ann. Probab.} {\bf 8} (1980), 539--575.

\bibitem{Sk}
A. V. Skorohod: {\em Random Processes with Independent Increments}.
Kluwer, Dordrecht, 1991.

\bibitem{So} R. Song: Estimates on the Dirichlet heat kernel
of domains above the graphs of bounded $C^{1,1}$ functions.
{\it Glasnik Math. \bf 39} (2004), 275-288.

\bibitem{SV06} R. Song and Z. Vondra\v{c}ek: Potential theory of special subordinators
and subordinate killed stable processes.  {\em J. Theoret. Probab.}
{\bf 19} (2006), 817--847.



\bibitem{SV07} R. Song and Z. Vondra\v{c}ek: Parabolic Harnack inequality
for the mixture of Brownian motion and stable process.
  \ {\em Tohoku Math. J.}
  {\bf 59} (2007), 1--19.

\end{thebibliography}
\end{document}